\documentclass[10pt]{article}
\usepackage[french,english]{babel}
\usepackage{amsmath}
\usepackage{amssymb}
\def\R{\mathbb{R}}
\def\N{\mathbb{N}}
\def\P{\mathbb{P}}
\def\E{\mathbb{E}} %
\def\L{\mathbb{L}}

\def\R{\mathbb{R}}

\def\Z{\mathbb{Z}}
\def\1{\mathbb{1}}
\def\1{\mbox{I\hspace{-.6em}1}} 
\def\var{\mbox{Var}\,}
\def\v{\mbox{Var}\,}

\def\cov{\mbox{Cov}\,}\def\Cov{\mbox{Cov}\,}
\def\lip{\mbox{Lip}\,}
\def\Lip{\mbox{Lip}\,}

\def\limiten{\renewcommand{\arraystretch}{0.5}
\begin{array}[t]{c}\stackrel{}{\longrightarrow} \\
{\scriptstyle n\rightarrow \infty} \end{array}
\renewcommand{\arraystretch}{1}}

\def\limiteloi{\renewcommand{\arraystretch}{0.5}
\begin{array}[t]{c}\stackrel{{\cal D}}{\longrightarrow} \\
{\scriptstyle n\rightarrow \infty} \end{array}
\renewcommand{\arraystretch}{1}}

\def\limiteprobaN{\renewcommand{\arraystretch}{0.5}
\begin{array}[t]{c}\stackrel{ \P}{\longrightarrow} \\
{\scriptstyle N\rightarrow \infty} \end{array}
\renewcommand{\arraystretch}{1}}

\def\limiteL{\renewcommand{\arraystretch}{0.5}
\begin{array}[t]{c}\stackrel{}{\longrightarrow} \\
{\scriptstyle L\rightarrow \infty} \end{array}
\renewcommand{\arraystretch}{1}}

\def\limiteas{\renewcommand{\arraystretch}{0.5}
\begin{array}[t]{c}\stackrel{a.s.}{\longrightarrow} \\
{\scriptstyle n\rightarrow \infty} \end{array}
\renewcommand{\arraystretch}{1}}
\author{Jean-Marc Bardet$^a$\footnote{Author for correspondence.}, Paul Doukhan$^{ab}$, Jos\'e Rafael
Le\'on$^c\footnote{This author aknowledges the
program ECOS-NORD of Fonacit, Venezuela, for its support.}$\\
~\\
{\small $^a$ SAMOS-MATISSE-CES, Universit\'e Paris I, 90 rue de
Tolbiac, 75013
Paris, FRANCE}\\
\small {\tt bardet@univ-paris1.fr}\\
{\small $^{b}$ LS-CREST, Timbre J340, 3 avenue Pierre Larousse,
92240 Malakoff, FRANCE}\\
\small  {\tt doukhan@ensae.fr}
\\
{\small $^c$ Universidad Central de Venezuela, Escuela de Matem\'atica.}\\
{\small Ap. 47197 Los Chaguaramos, Caracas 1041-A VENEZUELA}\\
\small  {\tt jleon@euler.ciens.ucv.ve} }
\date{}
\title{Uniform limit theorems for the integrated periodogram of
weakly dependent time series and their applications to Whittle's
estimate}

\newcommand{\findem}{\hfill\hbox{\hskip 4pt
\vrule width 5pt height 6pt depth 1.5pt}\vspace{.5cm}\par}

\newtheorem{prop}{Proposition}
\newtheorem{cor}{Corollary}
\newtheorem{theo}{Theorem}
\newtheorem{lem}{Lemma}

\begin{document}
\maketitle


\begin{abstract}
We prove uniform convergence results for the integrated periodogram
of a weakly dependent time series, namely a law of large numbers and
a central limit theorem. These results are applied to Whittle's
parametric estimation. Under general weak-dependence assumptions we
derive  uniform limit theorems and asymptotic normality of Whittle's
estimate for a large class of models. For instance the causal
$\theta$-weak dependence property allows a new and unified proof of
those results for ARCH($\infty$) and bilinear processes. Non causal
$\eta$-weak dependence yields the same limit theorems for two-sided
linear (with dependent inputs) or Volterra processes.
\end{abstract}
{\small {\it Key words: Periodogram, Weak dependence, Whittle estimate.}}\
\\
\
\\
{\it Running title:} Uniform limit theorems for the integrated
periodogram and Whittle's estimate.
\section{Introduction}
Parametric estimation from a sample of a stationary time series is
an important statistic problem both for theoretical research and
for its practical applications to real data. Whittle's
approximation likelihood estimate is particularly attractive for
numerous models like ARMA, linear processes, etc. mainly for two
reasons: first, Whittle's contrast does not depend on the marginal
law of the time series but only on its spectral density, and
second, its computation time is smaller than other parametric
estimation methods such as exact likelihood. Numerous papers have
been written on this estimation method after Whittle's  seminal
paper and in particular Hannan (1973), Rosenblatt (1985), and
Giraitis and Robinson (2001), established, respectively, the
asymptotic normality for Gaussian and causal linear, strong mixing
and ARCH($\infty$) processes. The case of long range dependent
processes, which is not considered here, was studied by Fox and
Taqqu (1986), and Giraitis and Surgailis (1990) (note also that in
a semi-parametric frame, Robinson (1995) and recently Dalla {\it
et al.} (2006) proposed an estimator of the long-memory parameter
based on local Whittle estimation). The main goal of the present
paper is to provide a unified treatment of this asymptotic
normality for a very rich class of weakly dependent time
processes, including those previously mentioned, but also some
never studied non causal or non linear processes.\\
~\\
More precisely, let $X=(X_k)_{k\in\Z}$ be a zero mean fourth-order
stationary time series with real values. Denote by $(R(s))_s$ the
covariogram of $X$, and by $(\kappa_4(i,j,k))_{i,j,k}$ the fourth
cumulants of $X$ defined by
\begin{eqnarray*}
R(s)&=&\cov (X_0,X_s)=\E \, (X_0X_s),~~~\mbox{for}~s\in \Z,\\
\kappa_4(i,j,k)&=&\E X_0X_iX_jX_k -\E X_0X_i \E X_jX_k -\E X_0X_j \E
X_iX_k- \E X_0X_k \E X_iX_j,~~~\mbox{for}~(i,j,k) \in \Z^3.
\end{eqnarray*}
We will use the following assumption on $X$ \\ \\
{\bf Assumption M:} $X$ is such that
\begin{eqnarray}
\gamma=\sum_{\ell\in\Z}R(\ell)^2<\infty~~\mbox{and}~~\kappa_4=
\sum_{i,j,k}|\kappa_4(i,j,k)|<\infty.
\end{eqnarray}
The periodogram of $X$ is
$$
I_n(\lambda) =\frac1{2\pi \cdot n} \left|\sum_{k=1}^n
X_ke^{-ik\lambda}\right|^2,~~~\mbox{for}~\lambda \in [-\pi,\pi[.
$$

Now, let $g:\R \rightarrow \R$ a $2\pi$-periodic function such that
$g\in\L^2([-\pi,\pi[)$ and define
\begin{eqnarray*}
J_n(g)&=& \int_{-\pi}^\pi g(\lambda)I_n(\lambda)\,
d\lambda,~~~\mbox{the integrated periodogram of}~X
\\
\mbox{and}~~J(g)&=& \int_{-\pi}^\pi g(\lambda)f(\lambda)\, d\lambda,
\end{eqnarray*}
with $f$ denoting the spectral density of $X$ (that exists and is in
$\L^2([-\pi,\pi[)$ from Assumption M) defined by
\begin{eqnarray*}\label{f}
f(\lambda)=\frac 1 {2\pi} \sum _{k \in \Z} R(k)\,
e^{ik\lambda}~~\mbox{for}~\lambda \in [-\pi,\pi[.
\end{eqnarray*}
Recall that $\displaystyle{I_n(\lambda) = \frac 1 {2 \pi}
\sum_{|k|<n}  \widehat R_n(k)e^{-ik\lambda}~~ \mbox{with}~~\widehat
R_n(k)=\frac 1 n \sum_{j=1\vee (1-k)}^{(n-k)\wedge  n}X_jX_{j+k}}$,
which is a biased estimate of $R(k)$. Thus, the periodogram
$I_n(\lambda)$ could be a natural estimator of the spectral density;
unfortunately it is not a consistent estimator. However, once
integrated with respect to some $\L^2$ function, its behavior
becomes quite smoother and can allow an estimation of the spectral
density. A special case of the integrated periodogram is  Whittle's
contrast, defined as a function $\beta \to J_n(h_\beta)$, where
$h_\beta$ is included in a class of functions depending on the
vector of parameters $\beta$. Whittle estimator minimizes this
contrast. As a consequence, uniform limit theorems for the
integrated periodogram $J_n(\cdot)$ are the appropriate tools for
obtaining uniform limit theorems for Whittle's contrast, that imply,
under additional conditions concerning the regularity of the
spectral density, limit theorems
for Whittle's estimators.      \\ \\
A uniform strong law of large numbers of integrated periodograms on
a Sobolev-type space (included in the space of $2\pi$-periodic
$\L^2$-functions) is first established only under assumption M.
Additional assumptions on the dependence properties of the time
series have to be specified for establishing central limit theorems.
Our choice has been to consider time series satisfying weak
dependence properties introduced and developed in Doukhan and
Louhichi (1999). Numerous reasons may explain this choice. First,
this frame of dependence includes a lot of models like causal or non
causal linear, bilinear, strong mixing processes and also dynamical
systems. Second, they are independent of the marginal distribution
of the time series. Finally, they can be easily used in various
statistical contexts, in particular in the case of the
integrated periodogram which is a quadratic form.  \\
These uniform limit theorems can be compared with those obtained in
Dahlhaus (1988) or Mikosch and Norvaisa (1997). Roughly speaking,
the presented results are obtained under weaker
conditions on time series, but considering different functional spaces.  \\
~\\
Two frames of weak dependence are considered here. The first one
exploits a causal property of dependence, the $\theta$-weak
dependence property (see Dedecker and Doukhan, 2003). Under certain
conditions, the uniform limit theorems for integrated periodogram
and asymptotic normality of Whittle's estimate are established.
These general results are new and extend Hannan's (1973) and
Rosenblatt's (1985) classical results for causal linear or strong
mixing processes. For example, parametric and causal ARCH($\infty$)
or bilinear processes (a very general class of models introduced by
Giraitis and Surgailis (2002), see definition (\ref{bilinear})) are
considered; under certain conditions, the asymptotic normality of
Whittle's estimators for those two classes of models is established
with the same method. (The case of causal ARCH($\infty$), and
therefore of GARCH(p,q), was already treated by Giraitis and
Robinson, 2001, under less restrictive conditions. However, their
proof is {\it ad hoc}
and cannot be used in a more general context).~\\
~\\
The second type of dependence under consideration is $\eta$-weak
dependence. This property allows to derive central limit theorems
for non causal processes, see \cite{bdl1} or \cite{DW}. These
results can be applied for instance, to two-sided linear or Volterra
processes (see their definition below in Section \ref{Whittle}). Let
us remark that usual proofs of central limit theorems for the
integrated periodogram are established by considering increments of
martingales or asymptotic results for strong mixing processes, which
is not a method  adapted for non causal processes, even in the
simple case of two-sided linear models. The proof of our results is
a corollary of a general functional central limit theorem for
$\eta$-weakly dependent processes, established by using  Bernstein's
blocks method. Even if our results may be sub optimal in terms of
the conditions linking the moment assumption with the decay rate of
weak dependence of the time series, they however cover numerous
models and open new perspectives of treatment for non causal
processes.\\
~\\
The paper is organized as follows. In Section \ref{Utheorem},
uniform limit theorems are presented with some applications to time
series. Section \ref{Whittle} is devoted to limit theorems satisfied
by Whittle's estimators, that are applied to several examples of
causal and non causal processes. Section \ref{nonc} contains the
main proofs, and a useful lemma is presented in an appendix (Section
\ref{appb}).
\section{Uniform limit theorems} \label{Utheorem}
\subsection{Notations and assumptions}
Afterwards we shall use the zero mean random variables
$$
Y_{j,k}=X_jX_{j+k}-R(k),~~\mbox{for all}~~(j,k) \in \Z^2.
$$
We intend to work in a Sobolev space ${\cal H}_s$ of locally $\L^2$
and $2\pi-$periodic functions, defined by
\begin{eqnarray*}
{\cal H}_s=\{ g \in\L^2([-\pi,\pi[)\, /\, \|g\|_{{\cal
H}_s}<\infty\}~~\mbox{with}~~ \|g\|_{{\cal
H}_s}^2=\sum_{\ell\in\Z}(1+|\ell|)^{2s} \cdot
|g_\ell|^2,~~~\mbox{and}~s>1/2,
\end{eqnarray*}
where $g(\lambda)=\sum_{\ell\in\Z}g_\ell \, e^{i\ell \lambda}$. This
space ${\cal H}_s$ is included in the space $C^\star$ of continuous
and $2\pi-$periodic functions and $\|g\|_\infty =\sup_{[-\pi,\pi[}
|g| \le \sqrt c_s \cdot \|g\|_{{\cal H}_s}$ with
\begin{eqnarray}\label{c}
c_s=\sum_{\ell\in\Z} (1+|\ell|)^{-2s}.
\end{eqnarray}
As usual ${\cal H}_s'$ denotes the dual of ${\cal H}_s$ with norm
defined from the identity $\|T\|_{{\cal H}'_s}=\sup_{\|g\|_{{\cal
H}_s}\le1}|T(g)|$. Hence if $T\in{\cal H}'_s$
$$
\|T\|^2_{{\cal H}'_s}=\sup_{\|g\|_{{\cal
H}_s}\le1}|T(g)|^2=\sum_{\ell \in\Z}(1+|\ell|)^{-2s} \cdot
|T(e^{i\ell\lambda})|^2.
$$
We study the asymptotic behavior of $J_n-J$ uniformly in the
function space ${\cal H}_s$ or equivalently as elements of  the
Hilbert space ${\cal H}'_s$.
\subsection{Uniform Strong Law of Large Numbers for the integrated
periodogram} We develop a uniform strong law of large numbers
(Uniform SLLN) for the integrated periodogram $(J_n(g))_g$. An
important feature is that the assumptions are only stated here in
terms of cumulant sums. Thus we need no additional assumption on the
dependence of the sequence $X$.
\begin{theo}[ULLN]\label{LGN}
If $X$ satisfies Assumption M, then
\begin{eqnarray}\label{USLLN}
\|J_n-J\|_{{\cal H}'_s}\limiteas 0.
\end{eqnarray}
\end{theo}
\paragraph{Remark} Let $X$ be a $4$-th order causal linear process,
{\it i.e.} $X_n=\sum_{k=0}^\infty a_k \, \xi_{n-k}$ for $n \in \Z$,
where $(\xi_k)_{k\in \Z}$ is a sequence of zero mean i.i.d.r.v. such
that ${\sum_kk\, a_k^2<\infty}$, ${\sum_k| a_k|<\infty}$ and $\E
\xi_0^4<\infty$. From Rosenblatt (1985, p. 59), and the
multilinearity of cumulants ({\it i.e.}
$|\kappa_4|\le\sum_{i,j,k}|\kappa_4(i,j,k)| = |c_4|\big(\sum_k|
a_k|\big)^4$ where $c_4$ is the $4$-th order cumulant of
$(\xi_k)_{k\in \Z}$), then $X$ satisfies Assumption M and therefore
the previous uniform SLLN yields. In Mikosch and Norvaisa (1997),
uniform LLN for the integrated periodogram are also proved for
causal linear processes. The integrated periodogram is then
considered on a general class of functions ${\cal F}$ endowed with a
pseudometric space of ${\cal L}^2$. In addition of the fourth moment
condition, the assumptions are respectively in Theorems 3.3 and 3.4,
$\sum_{k=0}^\infty k \, a_k^2<\infty$ for a uniform weak LLN and
$\sum_{j=1}^\infty k^{3/2} \,a_k^2<\infty$ for a uniform strong LLN.
As a conclusion, the uniform strong LLN (\ref{USLLN}) is satisfied
by a larger class of causal linear processes but the functional
space ${\cal H}_s$ is different from ${\cal F}$ (with also a
different distance).
\subsection{Uniform Central Limit Theorem for the integrated
periodogram} Now, we would like to establish a uniform central limit
theorem (UCLT) for the integrated periodogram $J_n(\cdot)$ on the
space of functions ${\cal H}_s$. Assumption M is not sufficient for
such a result. The dependence between the terms of the time series
$X$ has to be specified, and we will consider $2$ cases. Before
this, under Assumption M, we define for any $\lambda,\mu,\nu\in\R$,
the bispectral density
\begin{eqnarray*} \label{f4}
f_4(\lambda,\mu,\nu)=\frac1{(2\pi)^3}
\sum_{h=-\infty}^\infty\sum_{k=-\infty}^\infty\sum_{\ell=-\infty}^\infty
\kappa_4(h,k,\ell)e^{i(h\lambda+k\mu+\ell\nu)}
\end{eqnarray*}
the matrix $\Sigma=(\sigma_{\ell_i,\ell_j})_{1\leq i,j \leq m}$,
where $\ell_i$ are distinct integer numbers, with
\begin{eqnarray}
\label{sigma} \sigma_{k,\ell}= \sum_{h\in\Z}\Big(R(h)R(h+\ell-k)
+R(h+\ell)R(h-k)+\kappa_4(h,k,h+\ell)\Big ),
\end{eqnarray}
and for $g_1$ and $g_2$ in ${\cal H}_s$, the limiting covariance
$\Gamma(g_1,g_2)$
\begin{eqnarray} \label{Gamma}
\Gamma(g_1,g_2)=  4\pi \int_{-\pi}^\pi
g_1(\lambda)g_2(\lambda)f^2(\lambda)\, d\lambda +2\pi
\int_{-\pi}^\pi \int_{-\pi}^\pi
g_1(\lambda)g_2(\mu)f_4(\lambda,-\mu,\mu)\, d\lambda\, d\mu.
\end{eqnarray}
\subsubsection{UCLT for causal time series}\label{secUCLT}
This first case follows a classical methodology: the UCLT results
from the finite dimensional convergence and the tightness of the
process $Z=(Z_n(g))_{g \in {\cal H}_s}$ where
$\displaystyle{Z_n(g)=\sqrt n \, (J_n(g)-J(g))}$ for $n\in\N^*$ and
$g \in {\cal H}_s$. Since $g\mapsto Z_n(g)$ is a linear functional,
the finite dimensional convergence is a consequence of the
multidimensional central limit theorem for empirical covariances. \\
~\\
In the sequel, for $\ell \in \Z$, we will denote by ${\cal
M}^{(\ell)}_0$ a $\sigma$-algebra such that
\begin{eqnarray*}\label{M0}
{\cal M}^{(\ell)}_0 \supset \sigma\left (Y_{k,\ell},~k \leq 0\right
)=\sigma\left (X_{k}X_{k+\ell},~k \leq 0\right ).
\end{eqnarray*}
Here $\sigma(W_i,~i\in I)$ represents the $\sigma$-algebra generated
by $(W_i)_{i\in I}$. An example of such $\sigma$-algebra ${\cal
M}_0^{(\ell)}$ is ${\cal M}_0^{(\ell)}=
\sigma\left(X_{k_0},~X_{k_0+\ell} \right)$ for a given $k_0 \in \Z$.
The most classical of such $\sigma$-algebra is defined, for any
integer $p$ such that $p \geq \ell$,
\begin{eqnarray}\label{B0}
{\cal M}_0^{(\ell)}={\cal B}_p=\sigma\left (X_j,~j \leq p \right ).
\end{eqnarray}
\begin{lem}\label{lem3}
Let $(\ell_1,\ldots,\ell_m) \in\Z ^m$ be arbitrary distinct
non-negative integers ($m\in \N^*$). Let $X$ satisfy Assumption M
and be such that
\begin{eqnarray} \label{Proj}
\sum_{k\ge 0} \Big | \E \,\big (  Y_{0,\ell_i}\, \E \, (Y_{k,\ell_i}
\, |\, {\cal M}^{(\ell _i)}_0) \,\big ) \Big |<\infty~~~\mbox{for
all}~~i\in \{1,\ldots,m\}.
\end{eqnarray}
Then, if $\Sigma=(\sigma_{\ell_i,\ell_j})_{1\leq i,j \leq m}$
defined in (\ref{sigma}) is a nonsingular matrix,
$\mbox{with}~\limiteloi~\mbox{the weak convergence}$,
\begin{eqnarray}\label{RN}
\Big (\sqrt n\, (\widehat R_n(\ell_i)-\E\widehat R_n(\ell_i)) \Big
)_{1 \leq i \leq m}~~ \limiteloi {\cal N}_m(0,\Sigma).
\end{eqnarray}
\end{lem}
\paragraph{Remark} Lemma \ref{lem3} is a generalization of a result of
Rosenblatt (1985, Theorem 3, p. 58) which asserts the CLT if for
each $i \in \{1,\ldots,m\}$, $\displaystyle{ \sum_{k=0}^{\infty}
\left ( \E \Big ( \E \, (Y_{k,\ell_i}\, |\, {\cal B}_0) \Big
)^2\right )^{1/2}<\infty}$, and this implies (\ref{Proj}). An
analogue CLT in Hall and Heyde (1980) (Theorem 5.4, page 136) does
not seem to be adapted to work out the forthcoming examples.
\begin{theo}\label{UCLT}
Under assumptions of Lemma \ref{lem3}, the Uniform Central Limit
Theorem (UCLT) holds
\begin{eqnarray}\label{Z}
Z_n=\sqrt n \, (J_n-J)\limiteloi  Z~~\mbox{in the space}~{\cal
H}'_s,
\end{eqnarray}
with $(Z(g))_{g\in{\cal H}_s}$ the zero mean Gaussian process with
covariance $\Gamma(g_1,g_2)$ defined in (\ref{Gamma}).
\end{theo}
\paragraph{Remark} Dahlhaus (1988) and Mikosch and Norvaisa (1997)
established also UCLT of the integrated periodogram for a general
class of multivariate processes with high-order spectra (including
certain Gaussian processes) and for causal linear processes,
respectively. The considered functional spaces require to satisfy
conditions close to the classical Dudley's conditions. The UCLT
theorem \ref{UCLT} needs the simple  assumption $s>\frac12$, besides
the conditions for finite dimensional convergence. Recall that
assumption $s>\frac12$ is also required for obtaining uniform
equicontinuity, however additional dependence conditions are always
needed (see, for example, Dedecker and Louhichi, 2002). The book
edited by Dehling {\it et al.} (2002) contains also two other
articles on empirical spectral processes, concerned respectively
with uniform results for Gaussian processes (Dahlhaus and Polonik,
2002) and long range dependence processes (Soulier, 2002). Tails of
the variables and measures of entropy need to be controlled  as this
is usual for empirical processes (see the paper by Dedecker and
Louhichi in the same monograph); we do not have such assumption to
the price of very specific Hilbert type classes. An important
feature is the fact that assumption $s>\frac12$ is the best possible
assumption for empirical CLT in order that the corresponding Dudley
entropy integral converges. Uniform theorems on Hilbert classes are
indeed usually weaker results as stated in several works of Suquet
(see for instance Oliveira and Suquet, 1998); using De Acosta
(1970)'s criterion, entails that essentially the same assumptions
that those ensuring fidi convergence
usually  implies functional uniform convergence in Hilbert spaces. \\
~\\
{\bf Examples of time series satisfying Theorem \ref{UCLT}.} This
theorem is first applied to $3$ classical examples of time series,
which extend the known multidimensional CLT for integrated
periodogram. Then, the UCLT is established for a very rich class of
causal time series, the $\theta$-weakly dependent processes.
\begin{enumerate}
\item [\bf 1.]  Causal linear processes: let $X$ be a linear and causal
time series such that $X_n=\sum_{k=0}^\infty a_k\, \xi_{n-k}$ for $n
\in \Z$, where $(\xi_k)_{k\in\Z}$ is a sequence of zero mean
independent identically distributed random variables such that $\E
\xi_0^4<\infty$. Assume also that the real sequence $(a_k)$
satisfied ${\sum_kk\, a_k^2<\infty}$ and ${\sum_k| a_k|<\infty}$
(these conditions are weaker than those of Mikosch and Norvaisa,
1997, which are ${\sum_k k^{3/2+\varepsilon}\, a_k^2<\infty}$ with
$\varepsilon>0$, but the considered function spaces in both the UCLT
are different). Thus $X$ satisfies assumption M and from Rosenblatt
(1985, p. 59), we have for all $\ell \in \N$,
$\displaystyle{\sum_{k=0}^{\infty} \|\E \, (Y_{k,\ell}\,|\,{\cal
B}_\ell)\|_2<\infty}$ and thus (\ref{Proj}) is also satisfied. Then
the UCLT (\ref{Z}) holds.
\item [\bf 2.] Gaussian processes: let the sequence $(X_n)_{n\in\Z}$ be a
zero mean stationary Gaussian process such that $\sum_k
R(k)^2<\infty$. Then $X$ satisfies Assumption M for all $\ell \in
\Z$ and $k\in\N$, and, with the $\sigma$-algebra defined by
(\ref{B0}) for $m=k+\ell$, $\displaystyle{\big |\E \, \big (
Y_{0,\ell} \E \,(Y_{k,\ell}\,|\,{\cal M}^{(\ell)}_0) \big )\big |
= \Big |\E \, (X_0X_\ell X_k X_{k+\ell}) -R(\ell)^2 \Big |}$. But
$$\E \, (X_0X_\ell X_k X_{k+\ell})
-R(\ell)^2=R(k)^2+R(k+\ell)R(k-\ell)$$ for a zero mean stationary
Gaussian process. Thus,
\[
\sum_{k \geq 0} \big |\E \, \big (  Y_{0,\ell} \E
\,(Y_{k,\ell}\,|\,{\cal M}^{(\ell)}_0) \big )\big | \leq \sum_{k \in
\Z} R(k)^2 + \left ( \sum_{k \in \Z}R^2({k+\ell}) \right
)^{1/2}\hspace{-4mm} \cdot \hspace{-0mm} \left ( \sum_{k \in
\Z}R^2({k-\ell}) \right )^{1/2},
\]
from the Cauchy-Schwarz inequality for $\ell^2$ sequences. Therefore
the UCLT (\ref{Z}) applies and the assumption ($\sum_k
R(k)^2<\infty$) seems to be sharp in this case.
\item [\bf 3.] Strong mixing processes: here, we
consider the probability space $(\Omega,{\cal T},\P)$.
\begin{cor}\label{melange}
Let $X=(X_n)_{n\in \Z}$ be a sequence of random variables on
$(\Omega,{\cal T},\P)$ satisfying Assumption M. Assume that $X$ is
a $\alpha'$-mixing process, {\it i.e.} $$\displaystyle{\alpha '_n
=\sup_{\ell\ge 0}\left \{ \alpha \Big (\sigma(X_n,
X_{n+\ell})\,,\, {\cal B}_0 \Big ) \right \}\limiten 0},$$ where
$\displaystyle{\alpha \Big ({\cal A}\,,\, {\cal B} \Big )= \sup_{
A \in {\cal A},~B\, \in {\cal B}} \Big | \P (A \cap B)-\P (A) \P
(B) \Big |~}$ for ${\cal A},{\cal B} \subset {\cal T}$. Moreover,
with $Q_{|X|}$ denoting the quantile function of $|X|$ and
$\displaystyle{\alpha '^{-1}(u)= \sum_{k \geq 0} \1 _{\,u \leq
\alpha'_k} }$, assume that $\displaystyle{\int_0^{1}\alpha
'^{-1}(u) Q^4_{|X_0|}(u) \, du <\infty}$. Then the UCLT (\ref{Z})
holds.
\end{cor}
\paragraph{Remark} We note that
$\displaystyle{\alpha '_n \le \alpha_n= \alpha \Big (\sigma(X_k,\,k
\geq n)\,,\, {\cal B}_0 \Big )}$ the standard mixing coefficient in
Rosenblatt (1985). However, no simple counter example seems to be
available. Therefore, if $X$ is a strongly $\alpha$-mixing process
satisfying Assumption M such that
$\displaystyle{\int_0^{1}\alpha^{-1}(u) Q^4_{|X_0|}(u) \, du
<\infty}$, the UCLT (\ref{Z}) also holds.
\item [\bf 4.] Causal
$\theta$-weakly dependent processes: Doukhan and Louhichi (1999)
introduced the class of $\theta$-weakly dependent this notion was
developed later on in Dedecker and Doukhan (2003). It includes
numerous kinds of causal times series, for instance the strong
mixing processes (see other examples in Section \ref{Whittle}).
First, for $h:\R^u\to\R$ an arbitrary function, with $u \in \N^*$,
denote
$$
\Lip h=\sup_{(y_1,\ldots,y_u)\ne(x_1,\ldots,x_u)}  \frac{\left|
\,h(y_1,\ldots,y_u)- h(x_1,\ldots,x_u)\, \right|}
{|y_1-x_1|+\cdots+|y_u-x_u|} .
$$
The time series $X=(X_n)_{n\in \Z}$ is said $\theta-$weakly
dependent when there exists a sequence $(\theta_r)_{r\in \N}$
converging to $0$ such that for all $r\in \N$, all function $f:\R^2
\to \R$ satisfying $\|f\|_\infty\le 1$, and all random variable
$Z\in{\cal B}_0$ such that $\|Z\|_\infty<1$,
\begin{eqnarray}\label{theta}
|\cov(f(X_{j_1},X_{j_2}),Z)|\le 2 \cdot \lip f\cdot
\theta_r~~~\mbox{for all }~j_1,j_2\ge r.
\end{eqnarray}
Let us denote $\|X\|_m=(\E|X|^m)^{1/m}$.
\begin{cor}\label{cor_weak}
Let $X=(X_n)_{n\in \Z}$ a $\theta-$weakly dependent process
satisfying Assumption M.  We also suppose that there exists $m>4$
such that $\displaystyle{\|X_0\|_m<\infty~~\mbox{and}~~
\sum_{k=0}^\infty\theta_k^{\frac{m-4}{m-1}}<\infty}$. Then the UCLT
(\ref{Z}) holds.
\end{cor}
\end{enumerate}
\subsubsection{UCLT for non-causal weakly dependent time series}
>From the seminal paper of Doukhan and Louhichi (1999), a second
class of weakly dependent processes can be considered. This class
includes also non causal time series. A process $X=(X_n)_{n\in \Z}$
with values in $\R^d$ is a so-called $\eta-$weakly dependent process
if there exists a sequence $(\eta_r)_{r\in \N}$ converging to $0$,
satisfying
\begin{equation} \label{eta_ine}
\left|{{\Cov}}\Big( g_1(X_{i_1},\ldots,X_{i_u}),
g_2(X_{j_1},\ldots,X_{j_v})\Big ) \right| \leq\Big(u   (\lip g_1)
\|g_2\|_\infty+v   (\lip g_2)  \|g_1\|_\infty \Big )  \eta_{r}
\end{equation}
for all~$\displaystyle{\left \{ \begin{array}{l} \bullet~(u,v)\in
\N^*\times \N^*;\\
\bullet~(i_1,\ldots,i_u)\in \Z^u~\mbox{and}~(j_1,\ldots,j_v)\in
\Z^v~\mbox{with}~ i_1\leq \cdots \leq i_u < i_u+r\leq
j_1\leq\cdots\leq
j_v \\
\bullet~\mbox{functions}~ g_1:\R^{ud} \to \R~\mbox{and}~g_2:\R^{vd}
\to
\R~\mbox{satisfying}\\
\hspace{5cm}\|g_1\|_\infty\le \infty,~\|g_2\|_\infty\le
\infty,~\Lip g_1 < \infty~\mbox{and}~\Lip g_2 < \infty; \\
\end{array} \right . }$~\\
~\\
As a particular case of the functional limit theorem presented in
Bardet {\it et al.} (2005) a UCLT for integrated periodogram can
also be established, and more precisely a convergence rate to the
Gaussian law
\begin{theo}\label{Vite}
Let $X=(X_n)_{n\in \Z}$  be $\eta$-weakly dependent process and such
that assumption M holds. Suppose also that
$$
\exists m>4,~\mbox{such
that}~~\|X_0\|_m<\infty~~\mbox{and}~~\eta_n={\cal
O}(n^{-\alpha})~~\mbox{with}~~ \alpha>\max\Big (3\,;\,\frac
{2m-1}{m-4}\Big ).
$$
Then the UCLT (\ref{Z}) holds. Moreover, for $\phi:\R \to \R$ a
${\cal C}^3(\R)$ function having bounded derivatives up to order
$3$, and for $g \in {\cal H}_s$
\begin{eqnarray*}\label{vitesse}
\left  | \E \left [ \phi\Big ( \sqrt{n}(J_n(g)-J(g))\Big )-\phi \big
(\gamma(g) \cdot N \big ) \right ] \right | \le C \cdot n^{- \frac t
{t+3}\,\big (\frac {\alpha (m-4)-2m+1}{2(m+1+\alpha \cdot m)} \big )
}
\end{eqnarray*}
where $C>0$, $\displaystyle{t= \Big (\big (2 \alpha \, \frac
{m-2}{m-1} -1 \big ) \wedge \big (s- \frac {1} 2\big )  \Big ) }$,
$N\sim{\cal N}(0,1)$ and $\gamma^2(g) =\Gamma(g,g)$ defined in
(\ref{Gamma}).
\end{theo}
\begin{cor}\label{Vite2}
Under the same assumptions as in Theorem \ref{Vite}, for $\ell \in
\Z$ and $\phi:\R \to \R$ a ${\cal C}^3(\R)$ function having bounded
derivatives up to order $3$,
\begin{eqnarray*}\label{vitesse2}
\left  | \E \left [ \phi\Big (
\sqrt{n}(\widehat{R}_n(\ell)-R(\ell))\Big )-\phi \big
(\sigma_\ell\cdot N \big ) \right ] \right | \le C \cdot n^{-\frac
{\alpha (m-4)-2m+1}{2(m+1+\alpha \cdot m)}},
\end{eqnarray*}
with $C>0$, $N\sim{\cal N}(0,1)$ and
$$\displaystyle{\sigma_\ell^2=4 \pi \int _{-\pi} ^\pi
\cos^2(\lambda \ell))f^2(\lambda)d\lambda + 2 \pi \int _{-\pi} ^\pi
\int _{-\pi} ^\pi \cos ( \lambda \ell)\cos (\mu
\ell)f_4(\lambda,-\mu,\mu)d\lambda d\mu }.$$
\end{cor}

\paragraph{Remark} The convergence rates in both  the functional
limit theorems in Theorem \ref{Vite} and the Corollary \ref{Vite2}
are obtained by using Bernstein's blocks method. Even if the
conditions are perhaps not optimal we derive convergence rates
$n^{-\lambda}$ in those theorems with $\lambda < 1/2$. This
$\lambda$ could be made as close to $1/2$ as one wants by taking
$s$ and $\alpha$ large enough.
\section{Applications to parametric estimation}\label{Whittle}
Now we will apply the previous results to finite parameters
estimates. Let $X=(X_n)_{n \in \Z}$ be a time series satisfying
Assumption M. We denote by $f$ the spectral density  of $X$. Define
\begin{eqnarray}\label{norma}
\sigma^2=\exp\left(\frac1{2\pi}\int_{-\pi}^\pi\log f(\lambda)\,
d\lambda\right).
\end{eqnarray}
Following Rosenblatt (1985) we shall assume that the one-step
prediction error variance satifies $2\pi\sigma^2>0$. Assume that $f$
belongs to the family of functions defined in the form
\begin{eqnarray}\label{densite_beta}
f(\lambda)=f_{(\beta,\sigma^2)}(\lambda)=\sigma^2 \cdot g_\beta
(\lambda)~~~\mbox{for all}~~\lambda \in [-\pi,\pi[,
\end{eqnarray}
the function $f$ thus depends on a finite number of unknown
parameters, a variance term $\sigma^2$ and a
$\displaystyle\R^p$-vector $\beta$, where $\beta =
(\beta^{(1)},\ldots,\beta^{(p)})$. The normalization condition
(\ref{norma}) implies:
\begin{eqnarray}\label{normag}
\int_{-\pi}^\pi\log g_\beta(\lambda)\, d\lambda =0.
\end{eqnarray}
Denote also $\sigma^{*}$ and $\beta^* =
(\beta^{(1)*},\ldots,\beta^{(p)*})$ the true value of $\sigma$ and
$\beta$. As a consequence, for all $\lambda \in [-\pi,\pi[$, we will
now denote $\sigma^{*2} g_{\beta^*}(\lambda)$ the spectral density
of $X$. We will also consider the following conditions
\begin{itemize}
\item [$\bullet$] Condition C1: the true values $\sigma^*$ and $\beta^*$ are
such that $\sigma^*>0$ and $\beta^*$ lies in a region ${\cal K}
\subset \R^p$ where ${\cal K}$ is an open and relatively compact
set.
\item [$\bullet$] Condition C2: if $\beta_1$,
$\beta_2$ are distinct elements of ${\cal K}$, the set $\{ \lambda
\in [-\pi,\pi[,~g_{\beta_1}(\lambda) \neq g_{\beta_2}(\lambda)\}$
has positive Lebesgue measure.
\item [$\bullet$] Condition C3: for all $\beta \in {\cal K}$,
the function $\displaystyle{ \lambda \mapsto
g^{-1}_{\beta}(\lambda)=\frac 1 {g_{\beta} (\lambda)}\in {\cal
H}_s}$ and $\displaystyle{ \sup_{\beta \in {\cal K}} \|
g^{-1}_{\beta} \|_{{\cal H}_s} < \infty}$ with $s>1/2$.
\item [$\bullet$] Condition C4: for all $\lambda \in [-\pi,\pi[$,
the function $\beta \mapsto \displaystyle{ g^{-1}_{\beta}(\lambda)}$
is continuous on $\cal K$.
\item [$\bullet$] Condition C5: for all $\lambda \in [-\pi,\pi[$,
the function $\beta  \mapsto \displaystyle{
g^{-1}_{\beta}(\lambda)}$ is twice continuously differentiable on
${\cal K}$ and $\beta  \mapsto \displaystyle{ \int _{-\pi}^\pi \log
(g_\beta(\lambda))\, d\lambda}$ can be differentiated twice under
the integral sign.
\item [$\bullet$] Condition C6: there exists $s >1/2$ such that for all
$\beta \in {\cal K}$ and $(i,j)\in
\{1,\ldots,p\}$, the functions $\displaystyle{\lambda \mapsto \left
( \frac {\partial g^{-1}_{\beta}}{\partial \beta^{(i)}}\right
)_{\hspace{-2mm} \beta} (\lambda)}$, $\displaystyle{\lambda \mapsto
\left ( \frac {\partial ^2g^{-1}_{\beta}}{\partial \beta^{(i)}
\partial \beta^{(j)}}\right )_{\hspace{-2mm} \beta} (\lambda)}$ belong to
${\cal H}_s$ and $\displaystyle{ \sup_{\beta \in {\cal K}}\Big \|
\frac {\partial ^2g^{-1}_{\beta}}{\partial \beta^{(i)}
\partial \beta^{(j)}}\Big  \|_{{\cal H}_s} < \infty}$.
\item [$\bullet$] Condition C7: for all $\beta \in {\cal K}$, the function
$\displaystyle{ \lambda \mapsto g_{\beta}(\lambda)}$ is
continuously differentiable on $[-\pi,\pi[$.\\
\end{itemize}
Let $(X_1,\ldots,X_n)$ be a sample from $X$. As usually, define
Whittle's maximum likelihood estimators of $\beta^*$ and
$\sigma^{*2}$ as:
\begin{eqnarray*}\label{estim}
\widehat{\beta}_n = \mbox{Argmin} _{\beta \in {\cal K}} \left \{
J_n(g^{-1}_{\beta}) \right \}=\mbox{Argmin} _{\beta \in {\cal
K}}\left \{  \int_{-\pi}^\pi \frac
{I_n(\lambda)}{g_\beta(\lambda)}\, d\lambda \right
\}~~~\mbox{and}~~~\widehat{\sigma}^2_n =\frac 1
{2\pi}J_n(g^{-1}_{\widehat{\beta}_n})
\end{eqnarray*}
($\widehat{\beta}_n$ is supposed to be uniquely defined). In the
following paragraphs, we will show the strong consistency of the
estimators $\widehat{\beta}_n$ and $\widehat{\sigma}^2_n$.
\subsection{Asymptotic properties of Whittle's parametric
estimators}
\begin{theo}\label{estLGN}
Let $X$ satisfy the assumptions of Theorem \ref{LGN}. Under
Conditions C1-4, then
$$
\widehat{\beta}_n \limiteas \beta^*~~\mbox{and}~~
\widehat{\sigma}^2_n \limiteas \sigma^{*2}.
$$
\end{theo}
{\it Proof.} From Theorem \ref{LGN} and Condition C3 with
probability $1$,
$$
\lim _{n \to \infty} J_n(g^{-1}_{\beta}) = J(g^{-1}_{\beta}),
$$
uniformly on $\beta \in \cal K$. From Conditions C2 and
normalization condition ({\ref{normag}),
$$ J(g^{-1}_{\beta})>2
\pi \sigma^{*2}=J(g^{-1}_{\beta^*})~~\mbox{for all}~~\beta \neq
\beta^*
$$
(see Lemma 2, in Hannan, 1973). Therefore (see the details of the
proof of Theorem 1 in Hannan, 1973),
$\displaystyle{\widehat{\beta}_n=\mbox{Argmin} _{\beta \in {\cal K}}
\left \{ J_n(g^{-1}_{\beta}) \right \}}$ converges a.s. to $\beta^*$
and $\displaystyle{\widehat{\sigma}^2_n=\frac 1 {2 \pi}
J_n(g^{-1}_{\widehat{\beta}_n})}$ converges to
$\sigma^{*2}=\displaystyle{\frac 1 {2 \pi}J(g^{-1}_{\beta^*})}$.
\mbox{\findem}
\paragraph{Remarks on the conditions C1-4} The C1 and C2 conditions are
usual and can be found for example in Rosenblatt (1985) for mixing
time series or in Fox and Taqqu (1986) for strongly dependent times
series. The condition C4 is weaker than the condition of
differentiability generally required. It thus necessitates the
unusual condition C3 related to the uniform limit theorems \ref{LGN}
and
\ref{UCLT}.  \\
\begin{theo}\label{estTLC}
Let $X$ satisfy  either the assumptions of Theorem \ref{UCLT} or
those of Theorem \ref{Vite}. Under Conditions C1-6 and if the
matrix $W^*=(w_{ij}^*)_{1\leq i,j \leq p}$, with
\begin{eqnarray*}\label{W}
w_{ij}^*= \int _{-\pi}^\pi g^{2}_{\beta^*} (\lambda)   \Big ( \frac
{\partial g^{-1}_\beta} {\partial
\beta^{(i)}}\Big)_{\hspace{-1.3mm}\beta^*}\hspace{-2mm}(\lambda)
\Big ( \frac {\partial g^{-1}_\beta} {\partial
\beta^{(j)}}\Big)_{\hspace{-1.3mm}\beta^*}\hspace{-2mm}(\lambda) \,
d\lambda
\end{eqnarray*}
is nonsingular, then
\begin{eqnarray}\label{TLCbeta}
\sqrt{n} (\widehat{\beta}_n - \beta^*)& \limiteloi & {\cal N}_p\Big
(0\, , \,  (\sigma^*)^{-4} (W^*)^{-1} Q^*  (W^*)^{-1} \Big),
\end{eqnarray}
with a matrix $Q^*=(q_{ij}^*)_{1\leq i,j \leq p}$ whose entries
are defined by
\begin{eqnarray*}\label{Q}
q_{ij}^*=2 \pi  \left (2 \sigma^{*4} w_{ij}^*+  \int _{-\pi}^\pi
\int _{-\pi}^\pi f_4(\lambda,\mu,-\mu)  \Big ( \frac {\partial
g^{-1}_\beta} {\partial
\beta^{(i)}}\Big)_{\hspace{-1.3mm}\beta^*}\hspace{-2mm}(\lambda)
\Big ( \frac {\partial g^{-1}_\beta} {\partial
\beta^{(j)}}\Big)_{\hspace{-1.3mm}\beta^*}\hspace{-2mm}(\mu)\,
d\lambda\, d\mu \right ).
\end{eqnarray*}
\end{theo}
~\\ {\it Proof.} Let $U_n(\beta)=J_n(g^{-1}_{\beta})$. From
Conditions 2 and 5, $\beta \mapsto U_n(\beta)$ exists and is twice
differentiable on ${\cal K}$. Denote $\displaystyle{\frac {\partial
} {\partial \beta}U_n(\beta)}$ the vector $\displaystyle{\left (
\frac {\partial } {\partial \beta^{(i)}} U_n(\beta) \right) _{1\leq
i \leq p}}$ and $\displaystyle{ \frac {\partial^2 } {\partial
\beta^2}U_n(\beta) }$ the $(p \times p)$ matrix $\displaystyle{\left
(  \frac {\partial^2 } {\partial \beta^{(i)}
\partial \beta^{(j)}}U_n(\beta)\right) _{1\leq i,j \leq p}}$.
According to the mean value theorem,
$$
\frac {\partial } {\partial \beta} U_n(\widehat{\beta}_n)= \frac
{\partial } {\partial \beta} U_n(\beta^*)+\frac {\partial^2 }
{\partial \beta^2}U_n(\overline{\beta}_n)
(\widehat{\beta}_n-\beta^*),
$$
where $\| \overline{\beta}_n- \beta^* \|_p \leq \|\widehat{\beta}_n
- \beta^* \|_p$ (with $\| . \|_p$ denoting the euclidian norm in
$\R^p$). Since $\widehat{\beta}_n$ minimizes $\beta \mapsto
U_n(\beta)$, it follows that
\begin{eqnarray}\label{Un}
\frac {\partial } {\partial \beta} U_n(\beta^*)= \left [ - \frac
{\partial^2 } {\partial \beta^2}U_n(\overline{\beta}_n)\right ]
(\widehat{\beta}_n-\beta^*).
\end{eqnarray}
But, from Theorem \ref{estLGN}, $\widehat{\beta}_n \limiteas
\beta^*$ and then $\overline{\beta}_n \limiteas \beta^*$.
Consequently, from Condition C5, C6 and Theorem~\ref{LGN} (Uniform
Law of Large Number),
$$
\frac {\partial^2 } {\partial \beta^2}U_n(\overline{\beta}_n)
\limiteas  \left ( \int _{-\pi} ^\pi \frac {\partial^2 } {\partial
\beta^{(i)} \partial \beta^{(j)}} g^{-1}_{\beta^*}(\lambda) \cdot
\sigma^{*2} g_{\beta^*}(\lambda)\, d\lambda \right ) _{1\leq i,j
\leq p}=\sigma^{*2} W^*,
$$
(see Lemma 3 of Fox and Taqqu, 1986). Moreover, from Theorem
\ref{UCLT} and Condition C5,
\begin{eqnarray*}
\sqrt {n} \Big ( \frac {\partial } {\partial \beta} U_n(\beta^*) -
\frac {\partial } {\partial \beta}J(g^{-1}_{\beta^*}) \Big ) &
\limiteloi & {\cal N}_p(0,Q^*), \\
\mbox{and thus}~~~\sqrt {n}  \frac {\partial } {\partial \beta}
U_n(\beta^*) & \limiteloi & {\cal N}_p(0,Q^*),
\end{eqnarray*}
because~~ $\displaystyle{\frac {\partial } {\partial
\beta}J(g^{-1}_{\beta^*}) =\int _{-\pi} ^\pi \Big ( \frac {\partial
g^{-1}_{\beta}(\lambda) } {\partial \beta} \Big )_{\beta^*} \cdot
\sigma^{*2} g_{\beta^*}(\lambda)\, d\lambda =\sigma^{*2} \frac
{\partial } {\partial \beta} \left ( \int _{-\pi} ^\pi \log (
g^{-1}_{\beta}(\lambda))\, d\lambda \right)_{\beta^*} =0}$~~ from C6
and (\ref{normag}). Therefore, if the matrix $W^*$ is nonsingular,
from (\ref{Un}),
$$
\sqrt {n}(\widehat{\beta}_n-\beta^*) \limiteloi
-(\sigma^*)^{-2}(W^*)^{-1} \cdot{\cal N}_p(0,Q^*),
$$
and this completes the proof of Theorem \ref{estTLC}.
\mbox{\findem} ~\\
\begin{theo}\label{estTLCsigma}
Let $X$ satisfy either the assumptions of Theorem \ref{UCLT} or
those of Theorem \ref{Vite}. Under Conditions C1-7, then
\begin{eqnarray}
\label{TLCsigma} \sqrt{n} (\widehat{\sigma}^2_n - \sigma^*)&
\limiteloi & {\cal N}\Big (0\, , \, 2\sigma^{*4}+2 \pi \int
_{-\pi}^\pi \int _{-\pi}^\pi f_4(\lambda,\mu,-\mu)
g^{-1}_{\beta^*}(\lambda) g^{-1}_{\beta^*}(\mu)\, d\lambda\, d\mu
\Big).
\end{eqnarray}
Moreover, $\sqrt{n} (\widehat{\sigma}^2_n - \sigma^*)$ and $\sqrt{n}
(\widehat{\beta}_n - \beta^*)$ are jointly
asymptotically normal with covariance given by\\

$\displaystyle\lim _{n \to \infty} \sqrt{n} \Big ( \cov
(\widehat{\sigma}^2_n\, , \, \widehat{\beta}_n^{(i)}) \Big )_{1\leq
i \leq p}$ $$\quad\quad\quad\qquad = \big ( \sigma^{*2} \cdot W^*
\big )^{-1} \cdot \Big ( 2 \pi \int _{-\pi}^\pi \int _{-\pi}^\pi
f_4(\lambda,\mu,-\mu) g^{-1}_{\beta^*}(\lambda) \left ( \frac
{\partial } {\partial \beta^{(i)}} g^{-1}_{\beta^*}(\mu) \right)\,
d\lambda\, d\mu \Big )_{1\leq i \leq p}.
$$
\end{theo}
~\\ {\it Proof.} The Taylor's formula implies that
$$
U_n(\beta^*)=U_n(\widehat{\beta}_n)+(\beta^*-\widehat{\beta}_n)'
\cdot \Big (\frac {\partial^2 } {\partial
\beta^2}U_n(\underline{\beta}_n) \Big )\cdot
(\beta^*-\widehat{\beta}_n),
$$
with probability $1$, and with $\| \underline{\beta}_n- \beta^*
\|_p< \|\widehat{\beta}_n - \beta^* \|_p$. From previous Theorem
\ref{estTLC}, it follows
$$
\sqrt{n} ( U_n(\beta^*)- \sigma^{*2})
=\sqrt{n}(U_n(\widehat{\beta}_n)-\sigma^{*2})+ {\cal O}_p(n^{-1/2}).
$$
Under condition C7, $\E \Big ( U_n(\beta^*)\Big)=\sigma^{*2} + {\cal
O}(\log n /n)$ (see for instance Rosenblatt, 1985, p. 78) and thus
$\sqrt{n} \Big ( U_n(\beta^*)- \E \Big ( U_n(\beta^*)\Big) \Big )
\limiteloi {\cal N}\Big (0,\Gamma(g_{\beta^*}^{-1} ,
g_{\beta^*}^{-1}) \Big )$ with $g_{\beta^*}^{-1} \in {\cal H}_s$.
Therefore,
$$
\sqrt{n} \Big (U_n(\widehat{\beta}_n)-\sigma^{*2} \Big ) \limiteloi
{\cal N}\Big (0,\Gamma(g_{\beta^*}^{-1} , g_{\beta^*}^{-1}) \Big ),
$$
which implies relation (\ref{TLCsigma}). The end of the proof
follows the same arguments as in Rosenblatt (1985). \mbox{\findem}
\subsection{Examples of Whittle's parametric estimates for different time
series}
~\\
\underline{\bf Causal GARCH and ARCH($\infty$) processes}\\
~\\
The famous and from now on classical GARCH($q',q$) model was
introduced by Engle (1982) and Bollerslev (1986) and is given by
relations
\begin{eqnarray}\label{garch}
X_k= \rho_k \cdot \xi_k ~~~\mbox{with}~~ \rho_k^2 = a_0 +
\sum_{j=1}^q a_j X^2_{k-j}  + \sum_{j=1}^{q'} c_j \rho^2_{k-j},
\end{eqnarray}
where $(q',q)\in \N^2$, $a_0>0$, $a_j\geq 0$ and $c_j \geq 0$ for $j
\in \N$ and $(\xi_k)_{k\in \Z}$ are i.i.d. random variables with
zero mean (for an excellent survey about ARCH modelling, see
Giraitis {\it et al.}, 2005). Under some additional conditions, the
GARCH model can be written as a particular case of ARCH($\infty$)
model (introduced in Robinson, 1991) that satisfies
\begin{eqnarray}\label{archinfini}
X_k= \rho_k \cdot \xi_k ~~~\mbox{with}~~ \rho_k^2 = b_0 +
\sum_{j=1}^\infty  b_j X^2_{k-j},
\end{eqnarray}
with a sequence $(b_j)_j$ depending on the family $(a_j)$ and
$(c_j)$. Different sufficient conditions can be given for obtaining
a $m$-order stationary solution to (\ref{garch}) or
(\ref{archinfini}). Notice that for both models (\ref{garch}) or
(\ref{archinfini}), the spectral density is a constant. As a
consequence, the idea of Whittle's estimation in the GARCH case (see
Bollerslev, 1986) is based on the ARMA representation satisfied by
$(X_k^2)_{k \in \Z}$. Indeed, if $(X_k)$ is a solution of
(\ref{garch}) or (\ref{archinfini}), then $(X_k^2)$ can be written
as a solution of a particular case of equation (\ref{bilinear}) of
bilinear models  below (see Giraitis {\it et al.}, 2005). More
precisely,
$$
X_k^2=\varepsilon_k \Big (\gamma \cdot b_0+ \gamma \sum_{j=1}^\infty
b_j X^2_{k-j} \Big )+ \lambda_1 \cdot b_0 + \lambda_1
\sum_{j=1}^\infty b_j X^2_{k-j}~~~\mbox{for}~~k\in \Z,
$$
with $\varepsilon_k=(\xi_k^2-\lambda_1)/\gamma$ for $k\in \Z$,
$\lambda_1=\E \xi_0^2$ and $\gamma^2=\var (\xi_0^2)$. Moreover, the
time series $(Y_k)_{k \in \Z}$ defined by $
\displaystyle{Y_k=X_k^2-\lambda_1 \cdot b_0 \cdot \big ( 1
-\lambda_1 \sum_{j=1}^\infty b_j \big)^{-1} ~~~\mbox{for}~~k\in
\Z,}$ satisfies the forthcoming equation (\ref{bilinear}) with
parameter $c_0=0$ (as in Proposition \ref{estim_bilin} below).
Hence, a sufficient condition for the stationarity of $(X^2_k)_{k\in
\Z}$ with $\| X_0^2 \|_m<\infty$ is
$$
(\| \varepsilon_0 \|_m+1) \cdot \sum_{j=1}^\infty |b_j| <1
~~~\Longleftrightarrow ~~~ \Big (\frac {\| \xi^2_0-\lambda_1 \|_m}
\gamma+1 \Big ) \cdot\sum_{j=1}^\infty |b_j| <1.
$$
Let $(X_k)_{k \in \Z}$ a stationary solution of (\ref{archinfini})
and taking $\beta=(\beta^{(1)},\ldots,\beta^{(p)})\in \R^p$ such
that $b_j=b_j(\beta)$ for $j\in \N$, and $\sigma^2=\E
(X_0^2-\rho_0^2)=b_0^2(\beta) \cdot
h(\lambda_1,\gamma,\sum_{j=1}^\infty b_j(\beta))$, where $h$ is a
positive real function, the spectral density of $(X_k^2)_{k \in \Z}$
is
\begin{eqnarray}\label{densiteARCH}
f_{(\beta,\,\sigma^2)}(\lambda)=\frac{\sigma^2}{2 \pi} \cdot \Big |
1 -\sum_{j=1}^\infty b_j(\beta) \cdot e^{ij\lambda} \Big |^{-2}.
\end{eqnarray}
Then a Whittle's estimate can be used for estimating $\beta$ and
$\sigma^2$ parameters
in ARCH($\infty$) process: \\
\begin{prop}\label{estim_garch}
Let $X$ be a stationary ARCH($\infty$) time series following
equation (\ref{archinfini}), such that there exists $m>8$
satisfying $\E (|\xi_0|^m)<\infty$, satisfying the stationarity
condition,
$$
\Big ( \big (\frac {\| \xi^2_0-\lambda_1 \|_{m/2}} {\|
\xi^2_0-\lambda_1 \|_2}+1 \big )~\wedge~\| \xi_0 \|_m^2 \Big )
\cdot\sum_{j=1}^\infty |b_j(\beta)| <1,~~\mbox{and~one of the two
following assumptions}
$$
\begin{itemize}
\item {\it Geometric decay}: $\forall j \in \N$, $0 \leq b_j(\beta)$ and
$\exists \mu \in ]0,1[$ such that $b_j(\beta) ={\cal O}(\mu^{-j})$;
\item {\it Riemannian decay}: $\forall j \in \N$, $b_j(\beta) \geq 0$,
$\displaystyle{\exists \nu >\frac {2m-9}{m-8}}$ such that
$b_j(\beta)={\cal O}(j^{-\nu})$.
\end{itemize}
Then, under Conditions C1-7, the central limit theorems
(\ref{TLCbeta}) and (\ref{TLCsigma}) holds.
\end{prop}
\begin{cor}
If there exists $m > 8$ such that $X$ is a $m$-order stationary
GARCH(q',\,q) time series satisfying equation (\ref{garch}), then
with $\beta=(a_1,\ldots,a_q,c_1,\ldots,c_{q'})$, the central limit
theorems (\ref{TLCbeta}) and (\ref{TLCsigma}) are satisfied.
\end{cor}
{\it Proof.} First, Doukhan {\it et al} (2005) have shown that
ARCH($\infty$) processes satisfy the $\theta$-weak dependence
property: in the ``Geometric decay'' case, $\theta_r={\cal O}(e^{-c
\sqrt{r}})$ with $c>0$ and in the ``Riemannian decay'' case, with
$\nu>2$, $\displaystyle{\theta_r={\cal O}\big ( r^{-\nu+1} \big )}$.
Now by applying the following Lemma \ref{thetah} (see
section~\ref{appb}) for $h(x)=x^2$ and thus $a=2$, we deduce that
$(X_k^2)_{k \in \Z}$ is a $\theta^{\frac {m-2}{m-1}}$-weak dependent
time series. Hence, by denoting $\theta'=(\theta'_k)_{k\in \Z}$ the
weak dependent sequence of $X^2$ it holds: in ``Geometric decay''
case $\theta'_r={\cal O}(e^{-c \sqrt{r}})$ with $c>0$ and in the
``Riemannian decay'' case, $\displaystyle{\theta'_r={\cal O}\big (
r^{-\frac{(\nu-1)(m-2)}{m-1}} \big )}$. Corollary \ref{cor_weak}
implies that 1/ in the "Geometric decay" case, for all $\mu$,
$(X_k^2)_{k \in \Z}$ satisfies the Uniform CLT (\ref{Z}), 2/ in the
"Riemannian decay" case, $(X_k^2)_{k \in \Z}$ satisfies the Uniform
CLT (\ref{Z}) if~~ $\displaystyle{ \frac{(\nu-1)(m-2)}{m-1} \cdot
\frac {m/2-4}{m/2-1}
>1}$, {\it i.e.} $\displaystyle{\nu
>\frac {2m-9}{m-8}}$. \\
Finally, for GARCH$(q',q)$ the spectral density is explicit
$\displaystyle{f_{(\beta,\sigma^2)}(\lambda)=\frac {\sigma^2}{2
\pi} \frac {|1-\Sigma_{j=1}^q c_j e^{ij\lambda}|^2}
{|1-\Sigma_{j=1}^q a_j e^{ij\lambda}-\Sigma_{j=1}^q c_j
e^{ij\lambda}|^2}}$ and satisfies conditions C1-7 (more precisely
conditions C3, C6 and C7 using ${\cal H}_s$ with $s>1/2$ are
checked). Then, the proof of the corollary is a special
case of the "Geometric decay" case.  \mbox{\findem} ~\\
\paragraph{Remark} Zaffaroni (2003) studied the question of the Whittle's
estimation of the parameter of a stationary solution of
(\ref{garch}), this result was improved later on by Giraitis and
Robinson (2001). In this paper, the conditions to assure the
asymptotic normality of the Whittle's estimate are better than our
Theorem (\ref{estim_garch}), in two senses: 1/ only the $m=8$ is
required; 2/ the  conditions on the sequence $(b_j(\beta))$ in the
general case of ARCH($\infty$) model are only $b_0(\beta)>0$ and
$b_j(\beta) \geq 0$ for $j \in \N^*$ and the stationarity condition
$\| \xi_0 \|_m^2 \cdot\sum_{j=1}^\infty |b_j(\beta)| <1$. However,
their method for establishing the central limit theorem for the
periodogram is essentially {\it ad hoc} and can not be used for non
causal or non linear time series. The recent book of Straumann
(2005) also provides an up-to-date and complete overview to this
question. Chapter 8 of this book exposed the results of Mikosch and
Straumann (2002) that considered the case of intermediate moment
conditions of order $>4$ and $<8$ only for special case of
GARCH(1,1) processes. The convergence rates are
slower than the present ones.  \\
~\\
\underline{\bf Causal Bilinear processes} \\
~\\
Now, assume that $X=(X_k)_{k\in \Z}$ is a bilinear process (see the
seminal paper of Giraitis and Surgailis, 2002) satisfying the
equation
\begin{eqnarray}\label{bilinear}
X_k=\xi_k \Big ( a_0 + \sum_{j=1}^\infty a_j X_{k-j} \Big ) +c_0+
\sum_{j=1}^\infty c_j X_{k-j}~~~\mbox{for}~k\in \Z,
\end{eqnarray}
where $(\xi_k)_{k\in \Z}$ are i.i.d. random variables with zero mean
and such that $\| \xi_0 \|_m < +\infty$ with $m\geq 1$, and $a_j$,
$c_j$, $j\in \N$ are real coefficients. Assume $c_0=0$ and define
the generating functions
$$
\begin{array}{ll}
A(z)=\sum_{j=1}^\infty a_j z^j & C(z)=\sum_{j=1}^\infty c_j z^j \\
G(z)=(1-C(z))^{-1}=\sum_{j=0}^\infty g_jz^j ~~~~&
H(z)=A(z)G(z)=\sum_{j=1}^\infty h_j z^j.
\end{array}
$$
If $\| \xi_0 \|_m \cdot \sum_{j=1}^\infty |h_j| <\infty$, for
instance when $\| \xi_0 \|_m \cdot \Big ( \sum_{j=1}^\infty |a_j| +
\sum_{j=1}^\infty |c_j|\Big ) <1$ (see Giraitis and Surgailis,
2002), there exists a unique zero mean stationary and ergodic
solution $X$ in $\L^m(\Omega,{\cal A}, \P)$ of equation
(\ref{bilinear}) (see Doukhan {\it et al.}, 2006). For $m\geq 2$,
the covariogram of $X$ is~$R(k)=a_0^2 \cdot \| \xi_0 \|_2 \cdot \Big
(1-\sum_{j=1}^\infty h_j^2\Big )^{-1} \, \sum_{j=0}^\infty g_j\,
g_{j+k}$ and satisfies $\sum_k |R(k)|<\infty$. If we assume that
there exists $\beta=(\beta^{(1)}, \ldots, \beta^{(p)})$ such that
for all $k \in \Z$, $a_k=a_k(\beta)$ and $c_k=c_k(\beta)$, the
spectral density of $X$ exists and satisfies
$$
f_{(\beta,\,\sigma^2)}(\lambda)=\frac {a_0^2(\beta) \cdot \sigma^2}
{2\pi \big (1-\sum_{j=1}^\infty h_j^2 (\beta)\big )} \,\sum
_{k=-\infty}^{\infty} \sum_{j=0}^\infty g_j(\beta)\, g_{j+k}(\beta)
\, e^{-ik\lambda},
$$
with $\sigma^2=\| \xi_0 \|_2^2$. Like in Doukhan {\it et al.}
(2006), we consider three different cases of the convergence rate to
zero of the sequences $(a_k)$ and $(c_k)$. Then, using the previous
results for $\theta$-weak dependent time series, Whittle's estimate
of parameters $\beta$ and $\sigma^2$ satisfies the following
proposition
\begin{prop}\label{estim_bilin}
Let $X$ be a stationary bilinear time series satisfying equation
(\ref{bilinear}) with $c_0=0$, $\E (|\xi_0|^m)<\infty$ with $m>4$
and such that ~$\| \xi_0 \|_m \cdot \Big ( \sum_{j=1}^\infty |a_j| +
\sum_{j=1}^\infty |c_j|\Big ) <1$. Moreover, assume that $X$
satisfies one of the $3$ following conditions
\begin{itemize}
\item {\it Finite case}: $\exists J \in \N$ such that $\forall j>J$,
$a_j(\beta)=c_j(\beta)=0$;
\item {\it Geometric decay}: $\exists \mu \in ]0,1[$ such that
$\sum _j |c_j(\beta)| \mu ^{-j} \leq 1$ and $\forall j \in \N$, $0
\leq a_j(\beta) \leq \mu^j$;
\item {\it Riemannian decay}: $\forall j \in \N$, $c_j(\beta) \geq
0$, and $\displaystyle{\exists \nu_1 >\frac {2m-5}{m-4}}$ such that
$a_j(\beta)={\cal O}(j^{-\nu_1})$ and $\exists \nu_2 >0$ such that
$\sum _j c_j(\beta)j^{1+\nu_2}<\infty$, with $\displaystyle{\left \{
\begin{array}{l}\displaystyle{ \nu_2
>\frac {(m-4)\delta}{(m-1)\delta-(m-4)\log 2 }} \\
\displaystyle{\delta=\log \Big (1 + \frac {1-\sum_j
|c_j(\beta)|}{\sum_j c_j(\beta) j^{1+\nu_2}} \Big )>\log 2 \, \frac
{(m-4)}{(m-1)} }\end{array} \right . }$.
\end{itemize}
Then, under Conditions C1-7, the central limit theorems
(\ref{TLCbeta}) and (\ref{TLCsigma}) holds.
\end{prop}
~\\
{\it Proof.} Doukhan {\it et al.} (2006) studied the three different
cases of the Proposition and deduced the $\theta$-weak dependence
behavior of $X$ in each case. In the "Finite" and the "Geometric
decay" cases, $\theta_r={\cal O}(e^{-c \sqrt{r}})$ with $c>0$, which
implies the conditions required in Corollary \ref{cor_weak}.
Therefore, under Conditions C1-7, the models satisfy the
central limit theorems (\ref{TLCbeta}) and (\ref{TLCsigma}). \\
In the "Riemannian decay" case, $\displaystyle{\theta_r={\cal O}\Big
( \big ( \frac r {\log r} \big )^d \Big )}$ with
$\displaystyle{d=\max \Big ( -(\nu_1-1)\,;\, -\frac {\nu_2 \cdot
\delta}{\delta+\nu_2 \cdot \log 2}\Big )}$. The CLT (\ref{Z}) holds,
from Corollary \ref{cor_weak}, whenever $\displaystyle{d\cdot \frac
{m-1}{m-4}<-1}$. Thus, if $\displaystyle{1-\nu_1< -\frac
{m-4}{m-1}}$ and $\displaystyle{-\frac {\nu_2 \cdot
\delta}{\delta+\nu_2 \cdot \log 2}< -\frac {m-4}{m-1}}$, {\it i.e.}
$\displaystyle{\nu_1
>\frac {2m-5}{m-4}}$ and $\displaystyle{\nu_2
>\frac {(m-4)\delta}{(m-1)\delta-(m-4)\log 2 }}$, the assumptions of
Corollary \ref{cor_weak} hold true, and under Conditions C1-7, the
central limit theorems (\ref{TLCbeta}) and (\ref{TLCsigma}) are also
satisfied. \mbox{\findem}
~\\ ~\\
\underline{\bf Non-causal (two-sided) linear processes} \\
~\\
Let $X$ be a zero mean stationary non causal (two-sided) linear time
series satisfying
\begin{eqnarray}\label{nclp}
X_k=\sum_{j=-\infty}^\infty a_j \xi_{k-j}~~~\mbox{for}~~k\in \Z,
\end{eqnarray}
with $(a_k)_{k \in \Z} \in \R^\Z$ and $(\xi_k)_{k \in \Z}$ a
sequence of zero mean i.i.d. random variables such that $\E
(\xi_0^2)=\sigma^2<\infty$ and $\E (|\xi_0|^m)<\infty$ with $m \geq
4$. We assume that there exists $\beta=(\beta^{(1)}, \ldots,
\beta^{(p)})$ such that for all $k \in \Z$, $a_k=a_k(\beta)$.
Moreover, we assume that $(a_k(\beta))_{k \in \Z}$ is such that
$a_k(\beta)={\cal O}(|k|^{-a})$ with $a>1$. Therefore the spectral
density of $X$ exists and satisfies
$$
f_{(\beta,\sigma^2)}(\lambda)=\frac{\sigma^2}{2 \pi} \left | \sum
_{k=-\infty}^{\infty} a_k(\beta) e^{-ik\lambda} \right |^2.
$$
Then the results of the previous paragraph concerning non causal
weak dependent processes can be applied.
\begin{prop}\label{estim_lin}
Let $X$ be a linear time series satisfying (\ref{nclp}) with
$(a_k(\beta))_{k \in \Z} \in \R^\Z$ and $(\xi_k)_{k \in \Z}$ a
sequence of i.i.d. random variables with zero mean such that $\E
(\xi_0^2)=\sigma^2$ and $\E (|\xi_0|^m)<\infty$ with $m > 4$. We
assume that $(a_k(\beta))$ is such that
$$
a_k(\beta)={\cal O}(|k|^{-a})~~~\mbox{with}~~~a>\max \Big \{ \frac 7
2 \, ; \, \frac {5m-6}{2(m-4)} \Big \}.
$$
Then, under Conditions C1-7, the central limit theorems
(\ref{TLCbeta}) and (\ref{TLCsigma}) are satisfied.
\end{prop}
{\it Proof.} A $\eta$-weak dependence condition for non causal
linear random fields could be found in Doukhan and Lang (2002, p.
3); under the previous assumptions, $X$ is a $\eta$-weak dependent
time series with the relation: $\displaystyle{\eta_{2r}^2= {\cal O}
\big(\sum_{|k|>r}a_k^2(\beta)\big ) ~\Longrightarrow~~\eta_{r}=
{\cal O}\big(\frac 1 {r^{a-1/2}} \big )}$. Proposition
\ref{estim_lin} is then a consequence of Theorem \ref{Vite}.
\mbox{\findem}
\paragraph{Remark} 1/ The Condition C7 of central limit theorem
(\ref{TLCsigma}) is automatically satisfied by the convergence rate
of $(a_k)$ and therefore is not required in Proposition
\ref{estim_lin}. \\
~\\
2/The especial case of ARMA process i.e.
$f_{(\beta,\sigma^2)}(\lambda)=|\frac{A(e^{i\lambda})}{B(e^{i\lambda})}|^2$,
with $A$ and $B$ polynomials
$$A(z)=\sum_{k=0}^qa_kz^k,\quad a_0\neq 0,\quad
B(z)=\sum_{k=0}^pb_kz^k,\quad b_0=1,$$ without zeros in the circle
and with no common factors, is interesting. We have
$\beta=(a_0,\ldots,a_q,b_0,\ldots,b_p)$ (Rosenblatt, 1985 page
206) and $\sigma^2$ could be computed by using Jessen formula's. A
factorization as in (\ref{densite_beta}) is then possible. In
this setting conditions C1-7 are readily checked.\\
~\\
3/ To our knowledge, the known results about asymptotic
behavior of Whittle's parametric estimators for non-gaussian
linear processes are essentially devoted to one-sided (causal)
linear processes (see for instance, Hannan, 1973, Hall and Heyde,
1980, Rosenblatt, 1985, Brockwell and Davis, 1988). In such a
case, the conditions on $(a_k)$ are Conditions C1-6, with: $\sum_k
k\, a_k^2<\infty$ for the UCLT and the existence of $\sum_k k\,
a_k e^{-ik\lambda}$ for Condition C7. It
holds true if $m=4$ and $a_k={\cal O}(|k|^{-a})$ with $a>2$. \\
~\\
4/ There exist very few results in the case of two-sided linear
processes. In Rosenblatt (2000, p. 52) a condition for strong
mixing property for two-sided linear processes was given, but some
restrictive conditions on the process were also required for
obtaining a central limit theorem for Whittle's estimators: the
distribution of random variables $\xi_k$ has to be absolutely
continuous with respect to the Lebesgue measure with a bounded
variation density, $m>4+2\delta$ with $\delta>0$ and the central
limit theorem is obtained by using a tapered periodogram, under
assumption $\sum_{m=1}^\infty
\alpha_{4,\,\infty}(m)^{\delta/(2+\delta)}<\infty$. Here
$\alpha_{4,\,\infty}(m)\ge \alpha_m$ denotes a strong mixing
coefficient defined now with four points in the future instead of
2 for $\alpha'_m$, the same remark following Corollary
\ref{melange} still holds. The case of strongly dependent
two-sided linear processes was also treated by Giraitis and
Surgailis (1990) or Horvath and Shao (1999), however with more
restrictive conditions than Conditions C1-6 and with $a_k={\cal
O}(|k|^{-a})$ for a fixed $-1<a<0$.
  \\
~\\
5/ In the case of causal linear processes, it is well known that
$$\displaystyle{\sqrt{n} (\widehat{\beta}_n - \beta^*) \limiteloi
{\cal N}_p\big (0\, , \,  2\pi \cdot (W^*)^{-1} \big)},$$
$\widehat{\sigma}^2_n$ is a consistent estimate of $\sigma^*$, and
$\displaystyle{\sqrt{n} (\widehat{\beta}_n - \beta^*) }$ and
$\displaystyle{\sqrt{n} (\widehat{\sigma}^2_n - \sigma^*)}$ are
asymptotically normal and independent.
~\\ ~\\
\noindent
\textbf{{\underline{Non-causal Volterra processes}}}\\
~\\
A zero mean and non causal process $X=(X_t)_{t\in \Z}$ is called a
non-causal Volterra process, if  it satisfies
\begin{eqnarray}\label{ncvp}
X_k=\sum_{p=1}^\infty Y_k^{(p)},~~~\mbox{with}~~~
Y_k^{(p)}=\sum_{\tiny \begin{array} {c}j_1<j_2<\cdots<j_p\\
j_1,\ldots,j_p\in \Z
\end{array}} a_{j_1,\ldots,
j_p}\xi_{k-j_1}\cdots\xi_{k-j_p},
\end{eqnarray}
where $(a_{j_1,\ldots, j_p}) \in \R$ for $p \in \N^*$ and
$(j_1,\ldots,j_p)\in \Z^p$, and $(\xi_k)_{k \in \Z}$ a sequence of
zero mean i.i.d. random variables such that $\E
(\xi_0^2)=\sigma^2<\infty$ and $\E (|\xi_0|^m)<\infty$ with $m > 4$.
Such a Volterra process is a natural extension of the previous case
of non-causal linear process. From Doukhan (2003), the existence of
$X$ and thus the stationarity in $\L^m$ relies on the assumption
$$
\sum_{p=0}^\infty  \sum_{\tiny \begin{array} {c}j_1<j_2<\cdots<j_p\\
j_1,\ldots,j_p\in \Z
\end{array}}
\left|a_{j_1,\ldots, j_p}\right|^m \|\xi_0\|_{m}^p<\infty.
$$
%
Assume that there exists $\beta=(\beta^{(1)}, \ldots, \beta^{(p)})$
such that for all $p \in \N^*$ and $(j_1,\ldots,j_p)\in \Z^p$,
$a_{j_1,\ldots,j_p}=a_{j_1,\ldots,j_p}(\beta)$. Then the spectral
density of $X$ exists and satisfies
\begin{eqnarray}
\nonumber f_{(\beta,\sigma^2)}(\lambda)=  \sum_{p=1}^{\infty}
\frac{(\sigma^{2})^p}{2 \pi}  \sum_{k=-\infty}^\infty
\sum_{\tiny \begin{array} {c}j_1<j_2<\cdots<j_p\\
j_1,\ldots,j_p\in \Z
\end{array}} a_{j_1,\ldots, j_p}(\beta) \cdot
a_{j_1+k,\ldots, j_p+k}(\beta) \cdot e^{-ik\lambda}
\end{eqnarray}
(this formula is provided by the computation of the covariances of
$X$; remark that the representation with strictly ordered indices
$j_1<j_2<\cdots<j_p$ is fundamental). Certain conditions on the
asymptotic behavior of the coefficients $a_{j_1,\ldots,j_p}(\beta)$
give the $\eta$-weak dependence property of $X$ and then the
asymptotic normality of estimators
$(\widehat{\beta}_n,\widehat{\sigma}^{2}_n)$
\begin{prop}\label{estim_volt}
Let $X$ be a non-causal zero mean stationary Volterra process
satisfying relation (\ref{ncvp}) where $(a^{(p)}_{j_1,\ldots, j_p})
\in \R$ for $p \in \N^*$ and $(j_1,\ldots,j_p)\in \Z^p$, and
$(\xi_k)_{k \in \Z}$ a sequence of zero mean i.i.d. random variables
such that $\E (\xi_0^2)=\sigma^2<\infty$ and $\E (|\xi_0|^m)<\infty$
with $m > 4$.

Moreover, assume that the process is in some finite order chaos
(i.e. $a_{j_1,\ldots,j_p}(\beta)=0$ for $p>p_0$) and
$a_{j_1,\ldots,j_p}(\beta)$ is such that
$$
a_{j_1,\ldots,j_p}(\beta)={\cal O} \Big (\max_{1\le i\le p} \{
|j_i|^{-a} \} \Big ) ~~~\mbox{with}~~~a>4+\max \Big \{ 0 \, ; \,
\frac {11-m}{m-4} \Big \}.
$$
Then, under Conditions C1-7, the central limit theorems
(\ref{TLCbeta}) and (\ref{TLCsigma}) are satisfied.
\end{prop}
\paragraph{Remark}
A sharper dependence assumption is $\eta_r={\cal O} \Big (r^{1-a}
\Big )$ where $a$ is submitted to the same restriction as before;
recall that
\begin{equation*}
\eta_r\le 2\sum_{k=1}^\infty
\sum_{\tiny \begin{array} {c}j_1<j_2<\cdots<j_k\\
j_1<-r/2, \mbox{or } j_k\ge r/2
\end{array}} \left|a^{(k)}_{j_1,\ldots,
j_k}\right|\|\xi_0\|_{1}^k<\infty.
\end{equation*}
~\\
{\it Proof.} From Doukhan (2003), under the previous assumptions,
$X$ is a weakly dependent process with
\begin{equation*}
\eta_r\le \sum_{p=1}^\infty
\sum_{\tiny \begin{array} {c}j_1<j_2<\cdots<j_p\\
j_1<-r/2, \mbox{or } j_p\ge r/2
\end{array}} \left|a^{(p)}_{j_1,\ldots,
j_p}\right|\|\xi_0\|_{1}^p<\infty  ~~~~\Longrightarrow~~\eta_{r}=
{\cal O}\big(\frac 1 {r^{a+1}} \big ).
\end{equation*}
Proposition \ref{estim_volt} is then a consequence of Theorem
\ref{Vite}. \mbox{\findem} \\
~\\
\underline{\bf Non-causal linear processes with dependent
innovations} \\
~\\
Let $X=(X_n)_{n\in \N}$ be a zero mean stationary non causal
(two-sided) linear time series satisfying equation (\ref{nclp}) with
a dependent innovation process. More precisely, let $(\xi_n)_{n\in
\Z}$ be a weakly dependent fourth order centered stationary process
verifying Assumption M and such that $\E (\xi_0^2)=\sigma^2<\infty$.
Assume that there exists $\beta=(\beta^{(1)}, \ldots, \beta^{(p)})$
such that for all $k \in \Z$, $a_k=a_k(\beta)$ with
$a_k(\beta)={\cal O}(|k|^{-a})$ and $a>1$. Denoting
$g_{(\beta,\sigma^2)}$ the spectral density of the process
$(\xi_n)_{n\in \Z}$, the spectral density of $X$ exists and
satisfies
\begin{eqnarray*}\label{nclpdsdep}
f_{(\beta,\sigma^2)}(\lambda)&=& g_{(\beta,\sigma^2)}(\lambda) \cdot
\left | \sum _{k=-\infty}^{\infty} a_k(\beta) e^{-ik\lambda} \right
|^2.
\end{eqnarray*}
For instance, the process $(\xi_n)_{n\in \Z}$ may be a causal or a
non-causal ARCH($\infty$) or bilinear process. Following the results
of Doukhan and Wintenberger (2006), if $(\xi_n)_{n\in \Z}$ is a
$\eta$-weakly dependent process, then $X$ is an $\eta$-weakly
dependent process (with a sequence $(\eta_r)_r$ that can be
deduced). Thus, the asymptotic normality of the Whittle estimate of
parameters $(\beta,\sigma^2)$ could be established
\begin{prop}\label{estim_linnocausaldep}
Let $X$ be a linear time series satisfying (\ref{nclp}) with
$(a_k(\beta))_{k \in \Z} \in \R^\Z$ and $(\xi_k)_{k \in \Z}$ a
$\eta^{(\xi)}$-weakly dependent process with zero mean, a spectral
density $g_{(\beta,\sigma^2)}$ depending only on parameters
$(\beta,\sigma^2)$, and such that $\E (\xi_0^2)=\sigma^2$ and $\E
(|\xi_0|^m)<\infty$. Moreover, we assume that
\begin{eqnarray*}
a_k(\beta)={\cal O}(|k|^{-a}),~~~~\eta^{(\xi)}_r={\cal
O}(r^{-b}),~~~\mbox{with}~~~b \cdot  \frac
{(a-2)(m-2)}{(a-1)(m-1)}>\max \Big \{ 3 \, ; \, \frac {2m-1}{m-4}
\Big \}.
\end{eqnarray*}
Then, under Conditions C1-7, the central limit theorems
(\ref{TLCbeta}) and (\ref{TLCsigma}) are satisfied.
\end{prop}
~\\ {\it Proof.}  Doukhan and Wintenberger (2006) proved, under
assumptions $a_k(\beta)={\cal O}(|k|^{-a})$ and
$\eta^{(\xi)}_r={\cal O}(r^{-b})$, that $X$ is a $\eta$-weakly
dependent process with $\displaystyle{\eta_r={\cal O}\Big
(r^{-b\cdot \frac {(a-2)(m-2)}{(a-1)(m-1)}}\Big )}$.
Proposition~\ref{estim_linnocausaldep} is then a consequence of
Theorem \ref{Vite}.~\mbox{\findem}
\section{Appendix: proofs}\label{nonc}
\subsection{Proof of Theorem \ref{LGN}}
\begin{lem}\label{lll} If $X$ satisfies Assumption M,
then
$$ n \cdot \max_{ \ell \ge
0} \left ( \v \Big(\widehat R_n(\ell) \Big )\right ) \le
\kappa_4+2\gamma .
$$
\end{lem}
{\it Proof of Lemma \ref{lll}.} To prove this result, we use the
identity
$$
\cov (Y_{0,\ell},Y_{j,\ell})=\kappa_4(\ell,j,j+\ell)
+R(j)^2+R(j+\ell)R(j-\ell),
$$
and deduce from the stationarity of $(Y_{j,\ell})_{j \in \Z}$ when
$\ell$ is a fixed integer
\begin{eqnarray*}
n \cdot \v \Big (\widehat  R_n(\ell)\Big )&\le& \frac 1 n
\sum_{j=1\vee (1-\ell)}^{(n-\ell)\wedge  n} \sum_{j'=1\vee
(1-\ell)}^{(n-\ell)\wedge n}|\cov
(Y_{j,\ell},Y_{j',\ell})| \\
&\le&  \sum_{j \in \Z}|\cov (Y_{0,\ell},Y_{j,\ell})| \\
& \leq & \sum_{j \in \Z}
\left(|\kappa_4(\ell,j,j+\ell)|+2R(j)^2\right) \\
&\leq & \kappa_4+2\gamma ,
\end{eqnarray*}
by using Cauchy-Schwarz inequality for $\ell^2$-sequences. \mbox{\findem}\\
\begin{lem}\label{ll}
If $X$ satisfies Assumption M, then
$$
\E\|J_n-J\|_{{\cal H}'_s}^2\le  \frac 3 n \Big ( \gamma + c_s \cdot
(\kappa_4+2\gamma ) \Big),~~~\mbox{where}~c_s~\mbox{is defined in
}~(\ref{c}).
$$
\end{lem}
{\it Proof of Lemma \ref{ll}.} Let
$g(\lambda)=\sum_{\ell\in\Z}g_\ell \, e^{i \ell \lambda}\in{\cal
H}_s$. As in Doukhan and Le\'on (1989), we use the decomposition
\begin{eqnarray} \label{decompI_n}
J_n(g)-J(g)=-T_1(g)-T_2(g)+T_3(g)~~\mbox{with}~\left \{
\begin{array}{l}
\displaystyle{T_1(g)=\sum_{|\ell|\ge n}R(\ell) \, g_\ell,}\\
\displaystyle{T_2(g)=\frac1n\sum_{|\ell|<n}|\ell|\, R(\ell) \, g_\ell,}\\
\displaystyle{T_3(g)=\sum_{|\ell|<n} (\widehat R_n(\ell)-\E \widehat
R_n(\ell))\, g_\ell}
\end{array}
\right .
\end{eqnarray}
Remark that $T_3(g)=J_n(g)-\E J_n(g)$. Thus, we obtain the
inequality
$$
\E\|J_n-J\|_{{\cal H}'_s}^2\le 3(\|T_1\|_{{\cal
H}'_s}^2+\|T_2\|_{{\cal H}'_s}^2+\E \|T_3\|_{{\cal H}'_s}^2).
$$
Cauchy-Schwarz inequality yields
\begin{eqnarray*}
\|T_1\|_{{\cal H}'_s}^2&\le& \sum_{|\ell|\ge n}(1+|\ell|)^{-2s} \,
R(\ell)^2\le \frac 1 n \sum_{|\ell|\ge n}R(\ell)^2,\\
\|T_2\|_{{\cal H}'_s}^2&\le&\frac1{n^2}\sum_{|\ell|<n}|\ell|^2 \,
(1+|\ell|)^{-2s} \,  R(\ell)^2 \le \frac1{n}\sum_{|\ell|<n}
R(\ell)^2.
\end{eqnarray*}
Hence, $\displaystyle{\|T_1\|_{{\cal H}'_s}^2+ \|T_2\|_{{\cal
H}'_s}^2\le \frac \gamma{n}}$. Lemma \ref{lll} entails
\begin{eqnarray*}
\|T_3\|_{{\cal H}'_s}^2&\le&\sum_{|\ell|<n}(1+|\ell|)^{-2s} \,
(\widehat R_n(\ell)-\E \widehat R_n(\ell))^2,
\\
\E\|T_3\|_{{\cal H}'_s}^2&\le&\sum_{|\ell|<n} (1+|\ell|)^{-2s} \, \v
(\widehat R_n(\ell) )\le \frac 1 n \sum_{|\ell|<n} (1+|\ell|)^{-2s}
\, (\kappa_4+2\gamma)\le \frac{c_s(\kappa_4+2\gamma)}n,
\end{eqnarray*}
with $c_s$ defined in (\ref{c}). We combine those results to
deduce Lemma \ref{ll}.\findem ~\\
{\it Proof of Theorem \ref{LGN}.} We prove this strong law of large
numbers from a weak $\L^2$-LLN and Lemma \ref{ll}. The scheme of
proof is analogue to the one in the standard strong LLN in $\L^2$.
Set $t>0$. First, we know that for all random variables $X$ and $Y$,
we have $\P \left (X+Y \geq 2t \right ) \leq \P \left (X\geq t
\right )+ \P \left (Y\geq t \right )$. Thus
\begin{eqnarray}
\nonumber \P\left(\sup_{ n\ge N}\|J_n-J\|_{{\cal H}'_s}\ge 2
t\right)&\le& \sum_{k=[\sqrt N]}^\infty \P(\|J_{k^2}-J\|_{{\cal
H}'_s}\ge t)
\\
\nonumber &+& \sum_{k=[\sqrt N]}^\infty\P\left(\sup_{k^2\le
n<(k+1)^2} \|J_n-J_{k^2}\|_{{\cal H}'}\ge t\right)
\\
\label{ANBN}&\le& A_N+B_N.
\end{eqnarray}
>From Bienaym\'e-Tchebychev inequality, Lemma \ref{ll} implies that
\begin{eqnarray}\label{C1}
A_N\le  \frac {C_1} {t^2} \cdot\sum_{k\ge \sqrt N}\frac1{k^2},
\end{eqnarray}
with $C_1 \in \R_+$. Now set $\widetilde R_n(\ell) =\widehat
R_n(\ell)-\E\widehat R_n(\ell)$. The fluctuation term $B_N$ is more
involved and its bound is based on the same type of decomposition as
(\ref{decompI_n}), because for $k^2 <n$
$$
J_n(g)-J_{k^2}(g) =-T'_1(g)+T'_2(g)-T'_3(g),
$$
\begin{eqnarray*}
\mbox{with now}~~T'_1(g)&=&\sum_{k^2\le |\ell|<n}R(\ell) \, g_\ell,\\
T'_2(g)&=&\frac1n\sum_{k^2\le |\ell|<n}|\ell|\, R(\ell) \,
g_\ell,\\
\mbox{ and}~~ T'_3(g)&=&\sum_{|\ell|<k^2} \widetilde R_{k^2}(\ell)
\,  g_\ell - \sum_{|\ell|<n} \widetilde R_n(\ell)\, g_\ell.
\end{eqnarray*}
As previously, $\displaystyle{\|T'_1\|_{{\cal
H}'_s}^2+ \|T'_2\|_{{\cal H}'_s}^2\le \frac{\gamma }{k^2}}$. \\ \\
Set $\displaystyle{L_k=\max_{k^2\le n<(k+1)^2}
\|J_n-J_{k^2}\|_{{\cal H}'_s}}$ and
$\displaystyle{T^*_k=\max_{k^2\le n<(k+1)^2}\|T'_3\|_{{\cal
H}'_s}}$. Then,
$$
B_N\le \sum_{k\ge \sqrt N}b_k,\quad\mbox{with}\quad
b_k=\P\left(L_k\ge t \right)\le \frac{\E (L_k^2)}{t^2} .
$$
Now
$$
\E (L_k^2)\le 3(\|T'_1\|_{{\cal H}'_s}^2+ \|T'_2\|_{{\cal H}'_s}^2+
\E\|T^*_k\|_{{\cal H}'_s}^2) \le \frac{3\gamma }{k^2}+3 \cdot
\E\|T^*_k\|_{{\cal H}'_s}^2.
$$
Then, for $k^2\le n<(k+1)^2$ and $\ell \in \Z$,
\begin{eqnarray*}
\widetilde R_n(\ell)&=&\frac{k^2}n \widetilde R_{k^2}(\ell)+
\Delta_{\ell,n,k}\\
\Delta_{\ell,n,k}&=&\frac 1 {n} \sum_{h=k^2\wedge
(k^2-\ell))+1}^{n\wedge (n-\ell)}
(X_hX_{h+\ell}-R(\ell))=\frac1{n}\sum_{h=k^2\wedge
(k^2-\ell))+1}^{n\wedge (n-\ell)}Y_{h,\ell}.
\end{eqnarray*}
Remark that $\widetilde R_{k^2}(\ell)=0$ if $k^2\le |\ell|\le n$ and
thus $\widetilde R_{n}(\ell)=\Delta_{\ell,n,k}$ in such a case. Also
note that
\begin{eqnarray*}
\Delta^*_{\ell,k}&=&\max_{k^2\le n<(k+1)^2}|\Delta_{\ell,n,k}| \le
\frac1 {k^2}\sum_{h=(k^2\wedge
(k^2-\ell))+1}^{ (k^2+2k)\wedge ((k^2+2k)-\ell)}|Y_{h,\ell}|\\
\mbox{and thus}~~\E(\Delta^*_{\ell,k})^2&\le &\frac 1{k^4} \,(2k)^2
\cdot \max_{(h,\ell) \in \Z^2} \left (\E (
|Y_{h,\ell}|^2)\right)  \\
&\le & \frac{4}{k^2} \, \E(|X_{0}|^4).
\end{eqnarray*}
Write
\begin{eqnarray*}
T'_3(g)&=&\sum_{|\ell|<k^2} \widetilde R_{k^2}(\ell)
\left(1-\frac{k^2}n \right)g_\ell -
\sum_{|\ell|<n} \Delta_{\ell,n,k}\, g_\ell\\
|T^*_k(g)|&\le & \frac{2}{k}\sum_{|\ell|<k^2} |\widetilde
R_{k^2}(\ell) \, g_\ell| + \sum_{|\ell|<(k+1)^2}
\Delta^*_{\ell,k}\,|g_\ell|,
\end{eqnarray*}
and we thus deduce
\begin{eqnarray*}
\E\|T_k^*\|^2_{{\cal H}'_s}&\le&2c \cdot
\left(\frac{4}{k^2}\max_{\ell \in \Z}\left ( \v (\widehat
R_{k^2}(\ell))\right)+\max_{\ell \in \Z}\left (
\E(\Delta^*_{\ell,k})^2\right) \right) \le \frac{c\cdot A}{k^2}
\end{eqnarray*}
for a constant $A>0$ depending on $\E|X_0|^4,\kappa_4,$ and $\gamma$
only. Hence $b_k\le 3(\gamma+A\cdot c)/(k^2t^2)$ is a summable
series and, with $C_2>0$,
\begin{eqnarray}\label{C2}
B_N \leq  \frac {C_2}{t^2} \cdot \sum_{k\ge \sqrt N}\frac1{k^2}.
\end{eqnarray}
Then, (\ref{ANBN}), (\ref{C1}) and (\ref{C2}) imply
$\displaystyle{\sup_{n \geq N} \|J_n-J\|_{{\cal H}'_s} \limiteprobaN
0}$ what is equivalent to $\displaystyle{\|J_n-J\|_{{\cal H}'_s}
\limiteas 0}$.~\findem
\subsection{Proofs of the section \ref{secUCLT}}
First let us recall the following classical lemma (see a proof in
Rosenblatt (1985) \cite{Ro}, p. 58):
\begin{lem}\label{lem2} If $X$ satisfies Assumption M and $(\ell,k)\in\Z^2$
be arbitrary integers, then
\begin{eqnarray}
\nonumber n \cdot  \cov (&& \hspace{-0.7cm} \widehat R_n(k),\widehat
R_n(\ell))\limiten \sigma_{k,\ell}.
\end{eqnarray}
\end{lem}
{\it Proof of Lemma \ref{lem3}.} Under condition ({\ref{Proj}),
the projective criterion, introduced in Dedecker and Rio (2000),
{\em i.e.} $\displaystyle \sum_{k\ge 0}{\E \, \Big |
Y_{0,\ell_i}\, \E \, \Big (Y_{k,\ell_i} \, |\,
\sigma(\,Y_{j,\ell_i},~j\leq 0)\Big ) \, \Big|<\infty~~~\mbox{for
all}~i\in \{1,\ldots,m\}}$, holds. Therefore  each $\widehat
R_n(\ell_i)$ satisfies the central limit theorem. Now, by
considering a a linear combination of
$(Y_{j,\ell_1},\ldots,Y_{j,\ell_m})$, denoted $Z_j$, the
projective criterion is also satisfied by $(Z_j)_{j \in \Z}$
yielding the multidimensional central limit theorem (\ref{RN}).
\mbox{\findem} ~\\
~\\
{\it Proof of Theorem \ref{UCLT}.} We first prove the following
lemma:
\begin{lem}[Tightness] \label{lem1}
If $X$ satisfies Assumption M then the sequence of processes
$(Z_n)_{n\in\N^*}$ is tight in ${\cal H}'_s$. \\
\end{lem}
{\it Proof of Lemma \ref{lem1}.} Following de Acosta (1970), for
showing the tightness we only need to prove that the sequence is
flatly concentrated, this means that
$$
\E \Big (\|p_LZ_n\|_{{\cal H}'_s}^2 \Big ) \limiteL 0,
$$
where $p_L:{\cal H}'_s\to F'_L$ denotes the orthogonal projection on
the closed linear subspace $F'_L\subset {\cal H}'_s$ generated by
$(e_\ell)_{|\ell|\ge L}$ with $e_\ell(\lambda)=e^{i\ell\lambda}$
(also $F_L\subset {\cal H}_s$ denote the subspace generated by
$(e_\ell)_{|\ell|\ge L}$). Then, for $L >0$,
\begin{eqnarray}\label{Zn}
\|p_LZ_n\|_{{\cal H}'_s}=\sup_{\|g\|_{{\cal H}_s}<1,~g\in
F_L}|Z_n(g)|.
\end{eqnarray}
Thus, for $g =\sum_{|\ell| \geq L}g_\ell e_\ell \in F_L$ and
$\|g\|_{{\cal H}_s}<1$, using again the decomposition
(\ref{decompI_n}), we obtain
$$
|Z_n(g)|^2\le 3n \cdot (|T_1(g)|^2+|T_2(g)|^2+|T_3(g)|^2).
$$
First, from a Cauchy-Schwarz inequality, we have
\begin{eqnarray*}
n (|T_1(g)|^2+|T_2(g)|^2)&\le& n \left ( \sum_{|\ell| \geq n\vee L}
(1+|\ell|)^{-2s}R(\ell)^2   \sum_{|\ell|
\geq n\vee L} (1+|\ell|)^{2s}g_\ell^2  \right .\\
&&\hspace{1cm}+\left .~\1_{\{L <n\}}  \frac 1 {n^2} \sum_{L\le
|\ell|<n} |\ell|^2(1+|\ell|)^{-2s}R(\ell)^2  \sum_{L\le |\ell|<n}
(1+|\ell|)^{2s}g_\ell^2 \right ) \\
& \leq &  n  (1+|n\vee L|)^{-2s} \sum_{|\ell| \ge n\vee L}R(\ell)^2+
\1_{\{L <n\}}
\sum_{L\le |\ell|<n}|\ell| (1+|\ell|)^{-2s}  R(\ell)^2\\
&\le &\gamma  \Big ( L  (1+L)^{-2s}\,  \1_{\{L \geq n\}}+L
(1+L)^{-2s} \, \1_{\{L <n\}} \Big ).
\end{eqnarray*}
Thus, we obtain
\begin{eqnarray}\label{T1T2} \sup_{\|g\|_{{\cal
H}_s}<1,~g\in F_L} n   \left (|T_1(g)|^2+|T_2(g)|^2 \right )
\limiteL 0.
\end{eqnarray}
Also note that
\begin{eqnarray*}
\sqrt{n}\, T_3(g)&=&{\sqrt{n}}\sum_{|\ell|<n} (\widehat R_n(\ell)-\E
\widehat R_n(\ell))g_\ell\\
n \, |T_3(g)|^2&\le&{n}\sum_{L\le |\ell|<n}
(1+|\ell|)^{-2s}(\widehat
R_n(\ell)-\E \widehat R_n(\ell))^2\\
\E \Big ( \sup_{\|g\|_{{\cal H}_s}<1,~g\in F_L}\hspace{-3mm} n
|T_3(g)|^2 \Big )&\le&\sum_{L\le |\ell|<n}(1+|\ell|)^{-2s} \cdot
\sup_\ell \left(n\, \v (\widehat R_n(\ell))\right)\le
\sum_{|\ell|\ge L}(1+|\ell|)^{-2s} \, (\kappa_4+3\gamma) .
\end{eqnarray*}
Since $\sum_{\ell \in \Z} (1+|\ell|)^{-2s} <+\infty$, we deduce
$\displaystyle{\E \Big ( \sup_{\|g\|_{{\cal H}_s}<1,~g\in
F_L}\hspace{-3mm} n |T_3(g)|^2 \Big ) \limiteL 0}$.  With
(\ref{T1T2}) and (\ref{Zn}), the proof ends.~\mbox{\findem} \\
~\\
Now the proof of Theorem \ref{UCLT} can be achieved. Indeed, the
tightness allows to establish the functional central limit theorem.
Moreover, $J_n(g)-\E J_n(g)=\sum_{\ell \in \Z}g_\ell \, (\widehat
R_n(\ell)-\E \,(\widehat R_n(\ell))$, and therefore, from
(\ref{sigma}), the limiting covariance of the process
$(Z_n)_{n\ge1}$ is given by (\ref{Gamma}). More  details of the
finite dimensional convergence of the process
$(Z_n(g_1),\ldots,Z_n(g_k))$ can be found in Rosenblatt (1985,
Corollary 2, p. 61). \mbox{\findem} ~\\
~\\
{\it Proof of Corollary \ref{melange}.} Respectively from Rio's
inequality (1994) and from $X$'s stationarity,  we have, for all
$\ell,k \in \N$
\begin{eqnarray*}
\|Y_{0,\ell} \,\E \, (Y_{k+\ell,\ell}~|~{\cal B}_\ell)\|_1 &\le & 2
\int_0^{\alpha(\sigma(Y_{k+\ell,\ell}),{\cal B}_\ell) }
Q_{|Y_{0,\ell}|}(u)\,Q_{|Y_{k+\ell,\ell}|}(u) \, du \\
& \leq & 2 \int_0^{\alpha(\sigma(X_{\ell},X_{k+\ell}),{\cal B}_0) }
Q^2_{|Y_{0,\ell}|}(u) \, du.
\end{eqnarray*}
Therefore, for all $\ell, k \in \N$,
\begin{eqnarray*}
\|Y_{0,\ell} \,\E \, (Y_{k+\ell,\ell}~|~{\cal B}_\ell)\|_1 &\le & 2
\int_0^{\alpha '_{k}} Q^2_{|Y_{0,\ell}|}(u) \, du.
\end{eqnarray*}
Consequently, for all $\ell \in \N$
\begin{eqnarray*}
\sum_{k\geq 0}\|Y_{0,\ell} \,\E \, (Y_{k+\ell,\ell}~|~{\cal
B}_\ell)\|_1 &\le & 2 \int_0^{1}\left ( \sum_{k \geq 0} \1 _{\,u
\leq \alpha'_k} \right )
Q^2_{|Y_{0,\ell}|}(u) \,du, \\
&\le & 2 \int_0^{1} (\alpha '(u))^{-1} Q^2_{|Y_{0,\ell}|}(u) \,du.
\end{eqnarray*}
But Lemma 2.1 in Rio (2000) provides
\begin{eqnarray*}
\int_0^{1} (\alpha '(u))^{-1} Q^2_{|Y_{0,\ell}|}(u) du &\leq&
\int_0^{1} (\alpha '(u))^{-1} \Big
(Q_{|X_{0}|}(u)Q_{|X_{\ell}|}(u)+Q_{|R(\ell)|}(u)\Big)^2 \, du
\\
& \leq &  \int_0^{1} (\alpha '(u))^{-1}
(Q^2_{|X_{0}|}(u)+|R(\ell)|)^2 \, du,
\end{eqnarray*}
and therefore if $\displaystyle{\int_0^{1}(\alpha '(u))^{-1}
Q^4_{|X_0|}(u) \, du <\infty}$, then $\displaystyle{\sum_{k\geq
0}\|Y_{0,\ell} \,\E \, (Y_{k,\ell}~|~{\cal B}_\ell)\|_1 < \infty }$
for all $\ell \in \N$. \mbox{\findem}\\
~\\
{\it Proof of Corollary \ref{cor_weak}.} We truncate the variables
$X_j=f_M(X_j)+g_M(X_j)$ where, for $x\in \R$, we set $f_M(x)=(x
\wedge M)\vee(-M)$ (then $f_M \in [-M,M]$) and $g_M(x)=x-f_M(x)$
satisfies $|g_M(x)|\le|x| \, \1_{|x|\geq M}$.
Note that $\lip f_M=1$ and $\|f_M\|_\infty=M$. \\ \\
Then, $Y_{k,\ell}= \Big ( f_M(X_k)\,f_M(X_{k+\ell})-R'(\ell) \Big )
+U_{k,\ell,M}$
with: \\
\begin{itemize}
\item [$\bullet$] $\displaystyle{R'(\ell)=\E f_M(X_k)f_M(X_{k+\ell})  }$; \\
\item [$\bullet$]
$\displaystyle{U_{k,\ell,M}=g_M(X_k)f_M(X_{k+\ell})+f_M(X_k)g_M(X_{k+\ell})
+g_M(X_k)g_M(X_{k+\ell})+R'(\ell)-R(\ell).}$ \\
\end{itemize}
Therefore, $\E (\,U_{k,\ell,M})=0$ and, for $m$ such that
$\|X_0\|_m<\infty$, we derive
\begin{eqnarray}\label{U1}
\nonumber \| U_{k,\ell,M} \|_1\hspace{-2mm} &\leq& \hspace{-2mm}M \E
\,  |X_{k} |\cdot \1_{|X_{k}|
>M}  + M \E \, |X_{k+\ell}| \cdot \1 _{|X_{k+\ell}|
>M}  \\
\nonumber  \hspace{-2mm} && \hspace{3cm} + \E \, \Big((|X_{k}| \cdot
\1 _{|X_{k}|
>M })^2\Big ) +
\Big|\E \, \Big (X_0X_\ell-f_M(X_0)f_M(X_\ell)\Big) \Big |\\
\nonumber \hspace{-2mm} &\leq& \hspace{-2mm} 2M \|X_0 \|_m \Big ( \P
(X_0 >M)
\Big ) ^{1-1/m}+\|X_0 \|_m^2\Big ( \P (X_0 >M) \Big ) ^{1-2/m}\\
\nonumber \hspace{-2mm} && \hspace{3cm} +\E
\,|f_M(X_0)g_M(X_\ell)|+\E
\,|g_M(X_0)f_M(X_\ell)|+\E \,|g_M(X_0)g_M(X_\ell)| \\
\nonumber \hspace{-2mm} &\leq& \hspace{-2mm} 4M \|X_0 \|_m\Big
(\frac{\E \, | X_0|^m}{M^m} \Big ) ^{1-1/m} +2 \|X_0 \|_m^2\Big (
\frac{\E \, |
X_0|^m}{M^m} \Big ) ^{1-2/m} \\
\hspace{-2mm} &\leq& \hspace{-2mm}6\cdot  M^{2-m}\cdot  \|X_0
\|_m^m,
\end{eqnarray}
from H\"older and Markov inequalities. By the same procedure, we
also obtain
\begin{eqnarray}\label{U2}
\nonumber  \| U_{k,\ell,M} \|^2_2 &\leq& 6 \Big ( \E \,
g_M^2(X_k)f_M^2(X_{k+\ell}) +\E \, f_M^2(X_k)g_M^2(X_{k+\ell})+
\E \, g_M^2(X_k)g_M^2(X_{k+\ell}) \Big ) \\
& \leq  & 18 \cdot M^{4-m}\cdot \|X_0 \|_m^m.
\end{eqnarray}
Let $h_M$ be the function such that $h_M(x,y)=f_M(x)f_M(y)-R'(\ell)$
for all $(x,y) \in \R^2$. Note that $\|h_M\|_\infty \leq 2 M^2$ and
$\lip h_M=M$. Moreover, for all random variable $W$ in $\L
^1(\Omega,{\cal A},\P)$,
$$\displaystyle{\| \E \, (W \, | \, {\cal
B}_\ell) \|_1=\sup _{Z \in \L ^\infty(\Omega,{\cal B}_\ell,\P)~,~\|Z
\|_\infty \leq 1} \Big | \E \, (W \cdot Z ) \Big |}.
$$
Therefore,
\begin{eqnarray*}
\nonumber \Big \|h_M(X_0,X_\ell) \cdot  \E \, \Big (
h_M(X_k,X_{k+\ell}) \, | \, {\cal B}_\ell \Big ) \Big \|_1 & = &
\Big \| \E \, (h_M(X_0,X_\ell) \cdot h_M(X_k,X_{k+\ell}) \, |
\, {\cal B}_\ell)  \Big \|_1 \\
\nonumber & = &  \sup _{Z \in \L ^\infty(\Omega,{\cal B}_\ell,\P)~,~
\|Z \|_\infty \leq 1} \Big | \E \, \Big( Z\cdot h_M(X_0,X_\ell)
\cdot h_M(X_k,X_{k+\ell})
\Big) \Big | \\
& \leq & \sup _{Z' \in \L ^\infty(\Omega,{\cal B}_\ell,\P)~,~  \|Z'
\|_\infty \leq 2M^2} \Big | \cov \Big (h_M(X_k,X_{k+\ell})\, , \, Z'
\Big ) \Big |.
\end{eqnarray*}
Consequently, from the inequality (\ref{theta}) and the stationarity
of $X$, for all $k \geq 0$,
\begin{eqnarray}\label{theta1}
\Big \|h_M(X_0,X_\ell) \cdot  \E \, \Big ( h_M(X_k,X_{k+\ell}) \, |
\, {\cal B}_\ell \Big ) \Big \|_1  \leq 4\cdot  M^3 \cdot \theta _{k
-|\ell|}.
\end{eqnarray}
Thus,
\begin{eqnarray*}\label{theta2}
\nonumber \Big \|Y_{0,\ell} \cdot  \E\, \Big(Y_{k,\ell}\, |\, {\cal
B}_\ell\Big )\Big \|_1 &\leq & \Big \|U_{0,\ell,M} \cdot \E \, \Big
(U_{k,\ell,M}\, |\, {\cal B}_\ell \Big ) \Big \|_1 + \Big  \|
h_M(X_0,X_\ell)    \E \,
\Big (U_{k,\ell,M}\, |\, {\cal B}_\ell \Big )\Big  \|_1  \\
\nonumber && +\Big  \| U_{0,\ell,M} \cdot   \E \, \Big (h_M(X_k,X_
{\ell+k})\, |\, {\cal B}_\ell \Big ) \Big  \|_1+\Big
\|h_M(X_0,X_\ell) \cdot  \E \, \Big ( h_M(X_k,X_{k+\ell}) \, | \,
{\cal B}_\ell \Big ) \Big \|_1 \\
\nonumber & \leq & \Big \|U_{k,\ell,M}\Big \|_2^2 + 4M^2 \Big
\|U_{k,\ell,M}\Big \|_1 +\Big \|h_M(X_0,X_\ell) \cdot  \E \, \Big (
h_M(X_k,X_{k+\ell}) \, | \,
{\cal B}_\ell \Big ) \Big \|_1 \\
& \leq & 42 \cdot M^{4-m}\cdot \|X_0 \|_m^m +4 \cdot M^3\cdot
\theta _{k -|\ell|},
\end{eqnarray*}
from (\ref{U1}), (\ref{U2}) and  (\ref{theta1}). Set
$M=\theta_{k-|\ell|}^{-1/(m-1)}$, the previous calculations  prove
that, if there exists some $m>4$ such that
$$
\sum_{k=0}^\infty\theta_k^{1-3/(m-1)}<\infty,\qquad\mbox{and}\qquad
\|X_0\|_m<\infty,
$$
then the Uniform CLT (\ref{Z}) holds. \mbox{\findem}
\subsection{Proofs of Theorem \ref{Vite} and Corollary \ref{Vite2}}
{\it Proof of Theorem \ref{Vite}.} Set $g \in {\cal H}_s$ with
$g(\lambda)=\sum_{\ell \in \Z } g_\ell e^{i\lambda \ell}$, let $k
\in \N^*$ and define $g^{(k)}(\lambda)= \sum_{|\ell|<k} g_\ell
e^{i\lambda \ell}$. One obtains for all ${\cal C}^3(\R)$ function
$\phi$ with bounded derivatives up to order $3$
$$
\left |  \E \left [ \phi \Big ( \sqrt{n}(J_n(g)-J(g))\Big )-\phi
\Big (\sigma(g) \cdot N \Big ) \right  ]\right | \le
D_{1,n}^{(k)}+D_{2,n}^{(k)}+D_{3,n}^{(k)},
$$
with
\begin{eqnarray*}
D_{1,n}^{(k)} &=& \left |  \E \left [\phi \Big (
\sqrt{n}(J_n(g^{(k)})-J(g^{(k)}))\Big )-\phi \Big (\sigma(g^{(k)})
\cdot N
\Big ) \right  ]\right |  \\
D_{2,n}^{(k)} &=& \left |  \E \left [\phi \Big (\sigma(g^{(k)})
\cdot N \Big )
-\phi \Big (\sigma(g) \cdot N \Big ) \right  ]\right |  \\
D_{3,n}^{(k)} &=& \left |  \E \left [ \phi \Big (
\sqrt{n}(J_n(g^{(k)})-J(g^{(k)}))\Big )-\phi  \Big (
\sqrt{n}(J_n(g)-J(g))\Big )  \right  ]\right | .
\end{eqnarray*}
{\bf {\em Term $D_{1,n}^{(k)}$}:} For $i=1,\ldots,n$, set
$x_i=(X_{i+\ell})_{|\ell| <k}$ a stationary random vector in
$\R^{2k-1}$. The function
$$
h(x_i)=\sum_{|\ell|<k}g_\ell(X_iX_{i+\ell}-R(\ell))~~\mbox{for}~i=1,\ldots,n,
$$
satisfies the assumption H (namely $\E  h(x_0)=0$ and there exists
$a\ge 1$ and $A=A(d) \ge 1$ such as for all $u,v\in\R^d$, $ |h(u)|
\le A ( |u|^a\vee1)$ and $ |h(u)-h(v)|\le  A
\left(\left(|u|^{a-1}+|v|^{a-1}\right)\vee1\right) |u-v| $, see in
\cite{bdl1}). Here $a=2$ and $A=A(2k-1)=2k-1$. Define also
$$
S_n^{(k)}=\frac1{\sqrt
n}\sum_{i=1}^nh(x_i)=\sqrt{n}(J_n(g^{(k)})-J(g^{(k)})).
$$
By applying Theorem 1 of Bardet {\it et al.} (2005) to this function
$h$, one obtains with $C_1>0$ and $\displaystyle{\lambda= \frac
{\alpha(m-4)-2m+1}{2(m+1+\alpha \cdot m} }$,
\begin{eqnarray}\label{D1}
D_{1,n}^{(k)} \leq C_1 \cdot k^3 \cdot n ^{-\lambda},
\end{eqnarray}
{\bf {\em Term $D_{2,n}^{(k)}$}:} From inequality
$|a-b|\le|a^2-b^2|/a$ for  $a>0$, $b\ge0$, and from the mean value
theorem,
$$
D_{2,n}^{(k)} \leq \frac{\|\phi'\|_\infty}{ \sigma ^2(g)}\left|
\sigma ^2(g) -\sigma^2(g^{(k)}) \right |.
$$
But, from the expression (\ref{Gamma}), we deduce
\begin{multline*}
\left| \sigma ^2(g) -\sigma^2(g^{(k)}) \right | \leq \left| \frac 1
\pi \int _{-\pi} ^\pi
(g^2(\lambda)-(g^{(k)}(\lambda))^2)f^2(\lambda) \, d\lambda
\right . \\
\left.+  2 \pi \int _{-\pi} ^\pi \int _{-\pi} ^\pi
(g(\lambda)g(\mu)-g^{(k)}(\lambda)g^{(k)}(\mu))f_4(\lambda,-\mu,\mu)d\lambda
d\mu  \right |.
\end{multline*}
With $g \in {\cal H}_s$, we have
$$
\|g-g^{(k)}\|_\infty \leq  \Big ( \sum _{|\ell| \geq k}
(1+|\ell|)^{-2s} \Big )^{1/2}  \Big ( \sum _{|\ell| \geq k}
(1+|\ell|)^{2s} g_\ell^2 \Big )^{1/2}  \leq \Big ( \sum _{|\ell|
\geq k} (1+|\ell|)^{-2s} \Big )^{1/2}   \| g \| _{{\cal H}_s}.
$$
Consequently, with also $\|g+g^{(k)}\|_\infty \leq 2 \Big (
\sum_{\ell \in \Z} (1+|\ell|)^{-2s}  \Big )^{1/2} \cdot  \| g \|
_{{\cal H}_s}$, there exists $C_2>0$ such that
\begin{eqnarray}\label{D2}
D_{2,n}^{(k)}\leq C_2 \cdot k^{\frac {1} 2-s},
\end{eqnarray}
{\bf {\em Term $D_{3,n}^{(k)}$}:} First, from a Taylor expansion,
$$
D_{3,n}^{(k)} \leq \frac 1 2 \cdot  \| \phi'' \| _\infty  \cdot n
\cdot \E \Big (J_n(g-g^{(k)})-J(g-g^{(k)}) \Big )^2.
$$
With the same decomposition as in the proof of Lemma \ref{ll}, one
obtains
\begin{eqnarray*}
\E \Big (J_n(g-g^{(k)})-J(g-g^{(k)}) \Big )^2 &\leq & 3  \Big ( \Big
(  \sum _{|\ell| \geq n} R(\ell) g_\ell \Big )^2+ \Big (\frac 1 {n}
\sum _{k  \leq |\ell| < n}|\ell | R(\ell) g_\ell \Big )^2 +
\\
&& \hspace{3cm} + \Big \| \sum _{k \leq |\ell| < n}g_\ell \Big
(\widehat{R}_n(\ell)-\E (\widehat{R}_n(\ell)) \Big ) \Big \|_2^2
\Big ).
\end{eqnarray*}
First,
$$
\Big (  \sum _{|\ell| \geq n} R(\ell) g_\ell \Big )^2 \leq
(1+n)^{-2s} \cdot \sum _{|\ell| \geq n} R(\ell) ^2 \cdot \sum
_{|\ell| \geq n} (1+|\ell|)^{2s} g_\ell^2 \leq \frac 1 n \cdot
\|g\|^2_{{\cal H}_s} \cdot \sum _{|\ell| \geq n} R(\ell) ^2.
$$
Using the weak dependence of $(X_i)_i$ and with the same method as
in the proof of Lemma 2 in Bardet {\it et al.} (2005) adapted to the
function $h(x)=x$ (therefore with $a=1$),
$$
|R(\ell)| \leq c \cdot \eta _{|\ell|}^{\frac {m-2}{m-1}} \leq c
\cdot |\ell|^{-\alpha \, \frac {m-2}{m-1}} ,
$$
from the rate $\eta_{|\ell|}={\cal O} (|\ell|^{-\alpha} )$ with
$\alpha>3$. As a consequence,
$$
\sum _{|\ell| \geq n} R(\ell) ^2 \leq c \cdot n ^{1 -2 \alpha \,
\frac {m-2}{m-1}}~~\mbox{and}~~ \Big (  \sum _{|\ell| \geq n}
R(\ell) g_\ell \Big )^2\leq c \cdot n ^{ -2 \alpha \, \frac
{m-2}{m-1}}.
$$
In the same way,
$$
\left(\frac 1 {n} \sum _{k \leq |\ell| < n}|\ell | R(\ell) g_\ell
\right)^2 \leq \frac 1 n \cdot \|g\|^2_{{\cal H}_s} \cdot \sum _{k
\leq |\ell| < n} R(\ell) ^2  \leq  \frac c n \cdot k^{1 -2 \alpha \,
\frac {m-2}{m-1}}.
$$
Finally,
\begin{eqnarray*}
\Big \| \sum _{k \leq |\ell| < n}g_\ell \Big (\widehat{R}_n(\ell)-\E
(\widehat{R}_n(\ell)) \Big ) \Big \|_2^2 &\leq & \Big ( \sum _{k
\leq |\ell| < n} |g_\ell | \Big ( \v
(\widehat{R}_n(\ell))\Big )^{1/2} \Big )^2 \\
&\leq & \max_{\ell \in \Z} \Big (\v (\widehat{R}_n(\ell)) \Big )
\cdot  \|g\|^2_{{\cal H}_s} \cdot \sum _{k \leq |\ell| < n}
(1+|\ell|)^{-2s}\\
&\leq & \frac 1 n \cdot (\kappa_4+2 \gamma)\cdot  \|g\|^2_{{\cal
H}_s} \cdot \sum _{k \leq |\ell| < n} (1+|\ell|)^{-2s}~~ \mbox{from
Lemma \ref{lll}}.
\end{eqnarray*}
Summing up, there exists $C_3>0$ such that
\begin{eqnarray}\label{D3}
D_{3,n}^{(k)} \leq C_3\cdot ( k^{1 -2 \alpha \, \frac
{m-2}{m-1}}+k^{1-2s} ).
\end{eqnarray}
Now, with (\ref{D1}), (\ref{D2}) and (\ref{D3}), we deduce by
considering $\displaystyle{t= \Big (2 \alpha \, \frac {m-2}{m-1} -1
\Big ) \wedge (s-\frac {1} 2 })$ and selecting $k$ such that
$k^{t+3}=n^{\lambda}$, that there exists $C>0$ such that
$$
\left |  \E \left [ \phi \Big ( \sqrt{n}(J_n(g)-J(g))\Big )-\phi
\Big (\sigma(g) \cdot N \Big ) \right  ]\right | \le  C \cdot n ^{-
\frac t {t+3}\, \lambda}. \mbox{\findem}
$$
{\it Proof of Corollary \ref{Vite2}.} Set
$g(\lambda)=e^{i\ell\lambda}$ in Theorem \ref{Vite}. Since this
function belongs to each space ${\cal H}_s$, it follows that the
terms $D_{2,n}^{(k)}$ and $D_{3,n}^{(k)}$ both vanish and the result
follows from the bound (\ref{D1}). \mbox{\findem}
\section{Appendix: a useful lemma}\label{appb}
For a weakly dependent process, the following auxiliary lemma shows
that a function of this process is also a weakly dependent process
and moreover provides a relation between the two weak dependence
sequences. For an example of the use of such a result, see the
paragraph devoted to causal ARCH($\infty$) time series.
\begin{lem}\label{thetah}
Let $(X_i)_{i\in \Z}$ be a $\R^d$-valued and $\L^p$-stationary time
series with $p>1$, $h:\R^d\to \R$ be a function such that $|h(x)|
\leq c \|x\|^a$ and $|h(x)-h(y)| \leq c  \|x-y\|
(\|x\|^{a-1}+\|y\|^{a-1})$ for $(x,y) \in \left(\R^d\right)^2$,
$\|\cdot\|$ a norm on $\R^d$, with $0<c$ and $1\le a<p$. Let
$(Y_i)_{i\in \Z}$ be the stationary times series defined by
$Y_i=h(X_i)$ for $i\in \Z$. Then
\begin{itemize}
\item If $(X_i)_{i\in \Z}$ is $\theta$-weakly dependent time series,
then $(Y_i)_{i\in \Z}$ is a stationary $\theta^Y$-weakly dependent
time series, such that $\forall r \in \N$, $\theta^Y_r = C  \cdot
\theta_r^{\frac{p-a}{p-1}}$ with a constant $C>0$;
\item If $(X_i)_{i\in \Z}$ is $\eta$-weakly dependent time series, then
$(Y_i)_{i\in \Z}$ is a $\eta^Y$-weakly dependent time series, such
that $\forall r \in \N$, $\eta^Y_r = C \cdot
\eta_r^{\frac{p-a}{p-1}}$ with a constant $C>0$.
\end{itemize}
\end{lem}
~\\
{\it Proof.} We first assume that $d=1$. Let $f:\R^u \to \R$ and
$g:\R^v \to \R$ two real functions such that $\Lip f<\infty$,
$\|f\|_\infty \leq 1$, $\Lip g<\infty$, $\|g\|_\infty \leq 1$.
Denote $x^{(M)}=(x \wedge M)\vee(-M)$ for $x \in \R$. For simplicity
we first assume that $v=2$. Let $i_1,\ldots,i_u,j_1,\ldots,j_v \in
\Z^{u+v}$ such that $i_1,\ldots,i_u \geq r$ and $j_1,\ldots,j_v \leq
0$ and denote $x_{\bf i}=(X_{i_1},\ldots,X_{i_u})$ and $x_{\bf
j}=(X_{j_1},\ldots,X_{j_v})$. We then define functions $F:\R^u \to
\R$, $F^{(M)}:\R^u \to \R$ and $G:\R^v \to \R$, $G^{(M)}:\R^v \to
\R$ through the relations, $F(x_{\bf
i})=f(h(X_{i_1}),\ldots,h(X_{i_u}))$, $F^{(M)}(x_{\bf
i})=f(h(X^{(M)}_{i_1}),\ldots,h(X^{(M)}_{i_u}))$ and $G(x_{\bf
j})=g(h(X_{j_1}),\ldots,h(X_{j_v}))$, $G^{(M)}(x_{\bf
j})=g(h(X^{(M)}_{ j_1}),\ldots,h(X^{(M)}_{j _v }))$. Then
\begin{eqnarray*}
|\cov (F(x_{\bf i}),G(x_{\bf j}))| &\leq & | \cov (F(x_{\bf
i}),G(x_{\bf j})-G^{(M)}(x_{\bf j}))| + |\cov
(F(x_{\bf i}),G^{(M)}(x_{\bf j}))| \\
&\leq & 2\E |G(x_{\bf j})-G^{(M)}(x_{\bf j}))|  +2\E |F(x_{\bf
i})-F^{(M)}(x_{\bf i})| +  |\cov (F^{(M)}(x_{\bf i}),G^{(M)}(x_{\bf
j}))|.
\end{eqnarray*}
The last relation comes from $\|f\|_\infty \leq 1$. But we also have
\begin{eqnarray*}
\E |G(x_{\bf j})-G^{(M)}(x_{\bf j}))|
&\leq & v \cdot \Lip g \cdot \E |h(X_{0})-h(X^{(M)}_{0})| \\
&\leq & 2c  \cdot   v  \cdot \Lip g \cdot \E \big (|X_0|^a \cdot \1
_{|X_0|>M}
\big ) ~~~~~\mbox{(from the assumptions on}~h),\\
&\leq & 2c  \cdot v  \cdot  \Lip g  \cdot \|X_0\|_p \cdot M^{a-p}
~~~~~\mbox{(from Markov inequality)}.
\end{eqnarray*}
The same thing holds for $F$. Moreover, the functions $F^{(M)}$ and
$G^{(M)}$ satisfy $\Lip F^{(M)}= \Lip F^{(M)}= c \cdot M^{a-1}$,
with $c>0$, and $\| F^{(M)}\|_\infty \leq 1$,  $\| G^{(M)}\|_\infty
\leq 1$. Thus, from the definition of the weak dependence of $X$
\begin{eqnarray*}
\left|\cov \big (F^{(M)}(x_{\bf i}),G^{(M)}(x_{\bf j} )\big )\right|
&\leq &C  v  \Lip f  M^{a-1} \theta_r, \ \mbox{with $u=2$,
under condition }\theta;\\
&\leq &C  (v  \Lip f+u  \Lip g)   M^{a-1} \eta_r, \ \mbox{ under
condition }\eta.
\end{eqnarray*}
Finally, we obtain respectively
\begin{eqnarray*}
|\cov (F(x_{\bf i}),G(x_{\bf j}))| &\leq & C   v\Lip
f   \big ( M^{a-1} \theta_r + M^{a-p} \big )\\
&\leq & C   (v  \Lip f+u  \Lip g) \big ( M^{a-1} \eta_r + M^{a-p}
\big ).
\end{eqnarray*}
By the optimal choice of $M=\theta_r^{1/(1-p)}$, we get respectively
\begin{eqnarray*}
|\cov (f(Y_{i_1},\ldots,Y_{i_u}),g(Y_{i_1},\ldots,Y_{i_v}))| &\leq &
C v \Lip f
\theta_r^{\frac{p-a}{p-1}};\\
&\le& C  (v \Lip f+u  \Lip g) \eta_r^{\frac{p-a}{p-1}}.
\end{eqnarray*}
Changes for the vector valued case are simple, here truncations
should only be considered componentwize and yield the same result up
to a constant.~$ \mbox{\findem}$

\end{document}